\documentclass[10pt,a4paper]{article}

\usepackage[latin1]{inputenc}
\usepackage[T1]{fontenc}
\usepackage[english]{babel}
\usepackage{amsfonts}
\usepackage{euscript,amsthm}

\newtheorem{conj}{Conjecture}[section]
\newtheorem{theo}[conj]{Theorem}
\newtheorem{prop}{Proposition}[section]

\newtheorem{coro}[prop]{Corollary}
\newtheorem{lemm}{Lemma}[section]

\newtheorem{defi}{Definition}[section]
\newtheorem{nota}[defi]{Notation}

\begin{document}

% algèbres normées et anneaux usuels

\newcommand{\N}{\mathbb{N}}
\newcommand{\Z}{\mathbb Z}
\newcommand{\R}{\mathbb{R}}
\newcommand{\Q}{\mathbb{Q}}
\newcommand{\C}{\mathbb{C}}
\renewcommand{\H}{\mathbb{H}}
\renewcommand{\O}{\mathbb{O}}
\newcommand{\F}{\mathbb{F}}
\renewcommand{\S}{\mathbb{S}}

\renewcommand{\a}{{\cal A}}
\newcommand{\az}{\a_\Z}
\newcommand{\ak}{\a_k}

\newcommand{\rc}{\R_\C}
\newcommand{\cc}{\C_\C}
\newcommand{\hc}{\H_\C}
\newcommand{\oc}{\O_\C}

\newcommand{\rk}{\R_k}
\newcommand{\ck}{\C_k}
\newcommand{\hk}{\H_k}
\newcommand{\ok}{\O_k}

\newcommand{\rz}{\R_Z}
\newcommand{\cz}{\C_Z}
\newcommand{\hz}{\H_Z}
\newcommand{\oz}{\O_Z}

\newcommand{\RR}{\R_R}
\newcommand{\CR}{\C_R}
\newcommand{\HR}{\H_R}
\newcommand{\OR}{\O_R}

\newcommand{\re}{\mathtt{Re}}

% structure

\newcommand{\dem}{\noindent \underline {\bf D\'{e}monstration :} }
\newcommand{\pr}{\noindent \underline {\bf Proof :} }
\newcommand{\indic}{\noindent \underline {\bf Indication :} }
\newcommand{\rem}{\noindent \underline {\bf Remarque :} }
\newcommand{\rek}{\noindent \underline {\bf Remark :} }
\newcommand{\fin}{\begin{flushright} \vspace{-16pt}
 $\bullet$ \end{flushright}}
\newcommand{\lpara}{
\ \vspace{3pt}

\noindent}
\newcommand{\para}{
\

\

\noindent}
\newcommand{\Para}{
\

\

\

\noindent}

\newcommand{\sectionplus}[1]{\section{#1} \vspace{-5mm} \indent}
\newcommand{\subsectionplus}[1]{\subsection{#1} \vspace{-5mm} \indent}

% algèbre

\newcommand{\dual}{{\bf v}}
\newcommand{\com}{\mathtt{Com}}
\newcommand{\rg}{\mathtt{rg}}

\newcommand{\g}{\mathfrak g}
\newcommand{\h}{\mathfrak h}
\renewcommand{\u}{\mathfrak u}
\newcommand{\n}{\mathfrak n}
\newcommand{\plie}{\mathfrak p}
\newcommand{\q}{\mathfrak q}
\newcommand{\liesl}{\mathfrak {sl}}
\newcommand{\so}{\mathfrak {so}}

% diagrammes de Dynkin

% o-o-o-o

%#1=nb d'o horizontaux * 2 + 1.5
%#2=nb d'o horizontaux - 1

\newcommand{\dynkina}[2]{
\setlength{\unitlength}{2.5mm}
\begin{picture}(#1,3)
\put(0,0){$\circ$}
\multiput(2,0)(2,0){#2}{$\circ$}
\multiput(.73,.4)(2,0){#2}{\line(1,0){1.34}}
\end{picture}}

\newcommand{\dynkinap}[3]{
\setlength{\unitlength}{2.5mm}
\begin{picture}(#1,3)
\put(0,0){$\circ$}
\multiput(2,0)(2,0){#2}{$\circ$}
\multiput(.73,.4)(2,0){#2}{\line(1,0){1.34}}
\put(#3,0){$\bullet$}
\end{picture}}

\newcommand{\dynkinapp}[4]{
\setlength{\unitlength}{2.5mm}
\begin{picture}(#1,3)
\put(0,0){$\circ$}
\multiput(2,0)(2,0){#2}{$\circ$}
\multiput(.73,.4)(2,0){#2}{\line(1,0){1.34}}
\put(#3,0){$\bullet$}
\put(#4,0){$*$}
\end{picture}}

%*        
% >o-o-o-o
%o
%#1=nb d'o horizontaux * 2 + 1.5
%#2=nb d'o horizontaux - 1
\newcommand{\dynkindpspinp}[2]{
\setlength{\unitlength}{2.5mm}
\begin{picture}(#1,1.5)(.15,0)
\put(2,0){$\circ$}
\multiput(4,0)(2,0){#2}{$\circ$}
\multiput(2.73,.4)(2,0){#2}{\line(1,0){1.34}}
\put(-.15,1.17){$\bullet$}
\put(-.15,-1.17){$\circ$}
\put(.6,1.5){\line(5,-3){1.5}}
\put(.6,-.64){\line(5,3){1.5}}
\end{picture}
\vspace{.2cm}
}

\newcommand{\dynkinep}[3]{
\setlength{\unitlength}{2.5mm}
\begin{picture}(#1,1.5)(0,0)
\put(0,0){$\circ$}
\multiput(2,0)(2,0){#2}{$\circ$}
\multiput(.73,.4)(2,0){#2}{\line(1,0){1.34}}
\put(#3,0){$\bullet$}
\put(4,-2){$\circ$}
\put(4.42,-1.28){\line(0,1){1.36}}
\end{picture}
\vspace{.3cm}}

\newcommand{\dynkinepp}[4]{
\setlength{\unitlength}{2.5mm}
\begin{picture}(#1,1.5)(0,0)
\put(0,0){$\circ$}
\multiput(2,0)(2,0){#2}{$\circ$}
\multiput(.73,.4)(2,0){#2}{\line(1,0){1.34}}
\put(#3,0){$\bullet$}
\put(4,-2){$\circ$}
\put(4.42,-1.28){\line(0,1){1.36}}
\put(#4,0){$*$}
\end{picture}
\vspace{.2cm}}

\newcommand{\poidsesix}[6]{
\hspace{-.12cm}
\left [
\begin{array}{ccccc}
{} \hspace{-.2cm} #1 & {} \hspace{-.3cm} #2 & {} \hspace{-.3cm} #3 &
{} \hspace{-.3cm} #4 & {} \hspace{-.3cm} #5 \vspace{-.13cm}\\
\hspace{-.2cm} & \hspace{-.3cm} & {} \hspace{-.3cm} #6 &
{} \hspace{-.3cm} & {} \hspace{-.3cm}
\end{array}
\hspace{-.2cm}
\right ]      }

% géométrie algébrique

\newcommand{\spec}{\mathtt{Spec}}
\newcommand{\proj}{\mathtt{Proj}}
\newcommand{\sz}{{\spec\ \Z}}
\newcommand{\p}{{\mathbb P}}
\renewcommand{\P}{{\mathbb P}}
\newcommand{\A}{{\mathbb A}}
\newcommand{\pz}{\p_\Z}
\newcommand{\co}{{\cal O}}
\newcommand{\cf}{{\cal F}}
\newcommand{\cg}{{\cal G}}

%texte

\newcommand{\ssi}{si et seulement si }
\renewcommand{\iff}{if and only if }
\newcommand{\tr}{{}^t}
\newcommand{\trace}{\mbox{tr}}
\newcommand{\scal}[1]{\langle #1 \rangle}
\newcommand{\im}{\mathtt{Im}}

% diagrammes

\newcommand{\suiteexacte}[3]{#1 \rightarrow #2 \rightarrow #3}
\newcommand{\surmap}{\rightarrow \hspace{-.5cm} \rightarrow}
\newcommand{\limiteinverse}{\lim_\leftarrow}
\newcommand{\liste}{\

\begin{itemize}}
\newcommand{\codim}{\mbox{codim}}
\newcommand{\point}{^{ ^\bullet} \hspace{-.7mm}}
\newcommand{\X}{\mathfrak X}

\newcommand{\res}[2]{\vspace{.15cm} 

\noindent
{\bf #1 :} {\it #2} \vspace{.15cm} 

\noindent}

\newcommand{\fonction}[5]{
\begin{array}[t]{rrcll}
#1 & : & #2 & \rightarrow & #3 \\
   &   & #4 & \mapsto     & #5
\end{array}  }

\newcommand{\fonctionratsansnom}[4]{
\begin{array}[t]{ccl}
#1 & \dasharrow & #2 \\
#3 & \mapsto    & #4
\end{array}  }

\newcommand{\fonc}[3]{
#1 : #2 \mapsto #3  }

\newcommand{\directlim}[1]{
\lim_{\stackrel{\rightarrow}{#1}}    }

\newcommand{\inverselim}[1]{
\lim_{\stackrel{\rightarrow}{#1}}    }

\newcommand{\suitecourte}[3]{
0 \rightarrow #1 \rightarrow #2 \rightarrow #3 \rightarrow 0 }

\newcommand{\matdd}[4]{
\left (
\begin{array}{cc}
#1 & #2  \\
#3 & #4
\end{array}
\right )   }

\newcommand{\matddr}[4]{
\left (
\begin{array}{cc}
\hspace{-.2cm} #1 & \hspace{-.2cm}#2  \\
\hspace{-.2cm} #3 & \hspace{-.2cm}#4
\end{array}
\hspace{-.2cm} \right )   }

\newcommand{\mattt}[9]{
\left (
\begin{array}{ccc}
{} #1 & {} #2 & {} #3 \\
{} #4 & {} #5 & {} #6 \\
{} #7 & {} #8 & {} #9
\end{array}
\right )   }

\newcommand{\matttr}[9]{
\left (
\begin{array}{ccc}
{} \hspace{-.2cm} #1 & {} \hspace{-.2cm} #2 & {} \hspace{-.2cm} #3 \\
{} \hspace{-.2cm} #4 & {} \hspace{-.2cm} #5 & {} \hspace{-.2cm} #6 \\
{} \hspace{-.2cm} #7 & {} \hspace{-.2cm} #8 & {} \hspace{-.2cm} #9
\end{array}
\hspace{-.15cm}
\right )   }

\title{On Mukai flops for Scorza varieties}
\author{Pierre-Emmanuel Chaput\\
Pierre-Emmanuel.Chaput@math.univ-nantes.fr\\
Laboratoire de Mathématiques Jean Leray UMR 6629\\
2 rue de la Houssinière - BP 92208 - 44322 Nantes Cedex 3
}
\maketitle

\begin{center}
{\bf Abstract}
\end{center}

I give three descriptions of the Mukai flop of type $E_{6,I}$, one in
terms of Jordan algebras, one in terms of projective geometry over the
octonions, and one in terms of $\O$-blow-ups.
Each description shows that it is very similar to certain 
flops of type $A$. The Mukai flop of type $E_{6,II}$ is also described.

\begin{center}
{\bf Introduction}
\end{center}

In this article, I study a class of birational transformations called
``Mukai flops''. Let $G/P$ be a flag variety. Recall \cite{richardson}
that the natural map $T^* G/P \rightarrow \g^*$, where $\g$ is
the Lie algebra of $G$, has image
the closure of a single nilpotent orbit.

Sometimes, it happens that for two parabolic subgroups $P,Q \subset G$,
the images in $\g^*$ of $T^* G/P$ and $T^* G/Q$ are equal to the same
orbit closure $\overline \co$, and that moreover,
the above maps are birational isomorphisms. We therefore get a
birational map 
$$
\begin{array}{ccccc}
T^* G/P & & \dasharrow & & T^* G/Q,\\
& \searrow & & \swarrow\\
& & \overline \co
\end{array}
$$
called a Mukai flop. 

Since $T^* G/P$ is a symplectic variety,
nilpotent orbit closures provide a wide class of
examples of symplectic singularities and were studied also for this reason.
If $\overline \co$ is a nilpotent orbit closure,
then B. Fu showed that any symplectic resolution of $\overline \co$ is
given by a map $T^* G/P \rightarrow \overline \co$ \cite{fu}.
On the other hand,
in \cite{namikawa}, it is proved that any
Mukai flop can be described using fundamental ones, when $P$ (and $Q$)
is a maximal parabolic subgroup~: $G$ is then of type $A,D_{2n+1}$ or
$E_6$. In some sense, this provides a complete
understanding of the different symplectic resolutions of $\overline \co$
and the relations between them.
\lpara

In fact, the classical fundamental flops, when $G$ is of type $A_n$ or
$D_{2n+1}$, are easy to describe. The only items which are not very well
understood in this matter are the fundamental 
Mukai flops of type $E_6$, and the
purpose of this article is to fill this gap. 

Along with this
motivation in birational geometry, these
flops are key ingredients for the definition of generalized dual
varieties \cite{dual}
for a subvariety of the homogeneous space $G/P$, when $G$ is
of type $E_6$ and $P$ is the parabolic subgroup corresponding to the
root $\alpha_1$ or $\alpha_3$, with Bourbaki's notations \cite{bourbaki}.

For example, an easy consequence of theorem \ref{cotangent} is
theorem 2.1 in \cite{hermitian}, which
generalizes the fact that the dual variety of the smooth quadric
in $\p V$ defined by an invertible
symmetric map $f:V \rightarrow V^*$ is the quadric
in $\p V^*$ defined by $f^{-1}$, when the usual
projective space $\p V$ is replaced by any Scorza variety 
(see subsection \ref{def_scorza} for the definition of Scorza
varieties; for example, a grassmannian of
2-dimensional subspaces of an even-dimensional fixed space,
and $E_6/P_1,P_1$ the parabolic subgroup of the adjoint group of type
$E_6$ corresponding to the root $\alpha_1$, are Scorza varieties).

Finally, a third motivation is the study of the geometry of
exceptional homogeneous spaces. For example, subsection
\ref{para_tangent_severi} starts a study of the geometric properties
of $E_6/P_3$,
with a rather detailed description of its tangent bundle.

\lpara

There are two flops of type $E_6$, denoted $E_{6,I}$ (then $P$
corresponds to the root $\alpha_1$ and $Q$ to $\alpha_6$)
and $E_{6,II}$ ($P=P_2,Q=P_5$). 
I give three
descriptions of the flop $E_{6,I}$~: one via the geometry of the
corresponding flag variety, one using Jordan algebras, and one using a
new class of birational transformations that I call $\O$-blow-ups.

In fact, these three constructions work uniformly for $G/P$ any Scorza
variety. 
This gives for example a common description of the flop
$$\dynkinap{15}{6}{2} \longleftrightarrow \hspace{.6cm}
\dynkinap{15}{6}{10}$$
and the flop
$$\dynkinep{9}{4}{0} \hspace{.4cm} \longleftrightarrow \hspace{.5cm}
\dynkinep{9}{4}{8}.$$
This allows to understand better the latter.

\para

I now describe more precisely the contents of this article.
Let $k$ be a field and let $x \in \p^n_k$. Then a 
non-vanishing tangent vector 
$t \in T_xX$ defines a unique line $l$ with the following properties~:
\begin{itemize}
\item
$x \in l$\ ,
\item
$t \in T_xl$.
\end{itemize}
Moreover, the rational map $t \mapsto l$ is clearly the quotient map
$T_xX \simeq \A^n_k \dasharrow \p^{n-1}_k$, where $\p^{n-1}_k$ denotes the
variety of lines in $\p^n_k$ through $x$. Dually, we have a similar 
rational map
$T^*_xX \dasharrow {(\p^{n-1}_k)}^\dual$.

Section \ref{scorza} is devoted to proving the same kind of results
when the variety $\p^n_k$ is replaced by a Scorza variety (see subsection
\ref{def_scorza}), which after \cite{projective} is considered as
a projective space $\p^n_\a$ over a composition algebra $\a$, so that
when $\a=k$, we recover $\p^n_k$. So, in
this section, I show theorems \ref{tangent} and \ref{cotangent},
which have the following interpretation in terms of projective
geometry over $\a$~: given a generalized projective space $\p^n_\a$
and a point $x \in \p^n_\a$, there is a rational quadratic map 
$\overline{\nu^+_x} : T_x \p^n_\a \dasharrow \p^{n-1}_\a$, which maps
a tangent vector to the unique $\a$-line through it. Dually, there is
a similar map 
$\overline{\nu^-_x} : T^*_x \p^n_\a \dasharrow {(\p^{n-1}_\a)}^\dual$.

Propositions \ref{fano_tangent} and \ref{fano_cotangent} show that
polarizing $\nu_x^+$ (resp. $\nu_x^-$), 
one gets an isomorphism between the
variety of lines in $\p^n_\a$ through $x$ and the Fano variety of
maximal linear subspaces included in $\p^{n-1}_\a$ 
(resp. ${(\p^{n-1}_\a)}^\dual$). These two results don't have analogs when
$\a = k$.

\lpara

Note that this last ${(\p^{n-1}_\a)}^\dual$ is the projective space of
hyperplanes containing $x$; it is therefore included in 
${(\p^n_\a)}^\dual$.
The connection with Mukai flops is as follows~: assume that 
$G/P$ is the Scorza variety $\p^n_\a$. 
Let $x \in G/P$. We will see that there is a Mukai flop
$T^* \p^n_\a \dasharrow T^* {(\p^n_\a)}^\dual$.
The structure map
$T^* G/Q \rightarrow G/Q$ and this Mukai flop yield a composition 
$T^*_x \p^n_\a = T^*_x G/P \dasharrow T^* G/Q \rightarrow G/Q 
= {(\p^n_\a)}^\dual$.
Theorem \ref{cotangent} shows that this composition is the map
$\overline{\nu^-_x} : T^*_x \p^n_\a \dasharrow {(\p^{n-1}_\a)}^\dual
\subset {(\p^n_\a)}^\dual$.

Then, I show a general canonical
isomorphism of quotients of tangent spaces to
homogeneous spaces (theorem \ref{mukai}). As a particular case, this
theorem gives a way of computing a Mukai flop $T^* G/P \dasharrow
T^*G/Q$ once we know the composition $T^* G/P \dasharrow G/Q$.
I deduce a
description of the flop of type $E_{6,I}$ (proposition
\ref{flop_premier}), in terms of Jordan algebras.

In subsection \ref{eclate}, I give a maybe more geometric
description of the flop $E_{6,I}$. Recall that the minimal resolution
of the simplest Mukai flop
$T^* \p^n \dasharrow T^* {(\p^n)}^\dual$ is the blow-up of $T^* \p^n$
along the zero section. I show that the same result holds for the
$E_{6,I}$-flop $T^* \p^2_\O \dasharrow T^* {(\p^2_\O)}^\dual$, if one
replaces the usual notion of blow-up with an octonionic version of it
(theorem \ref{eclatement}).

Finally, concerning the Mukai flop of type $E_{6,II}$, I use the fact
that the homogeneous space $E_6/P_3$ can be realized as the space of
lines included in $E_6/P_1$.  Theorem 
\ref{flop_second} uses this model and the study of the tangent bundle 
$T(E_6/P_3)$ performed in subsection \ref{para_tangent_severi} to give also
a description of the Mukai flop of type $E_{6,II}$.

\lpara

Sections \ref{spin10} and \ref{sl5} study the restriction of the flops
to a cotangent space in the two cases when $G$
is of type $E_6$. They are of course $L$-equivariant rational maps, if
$L$ is a Levi factor of $P$, and happen to be quite subtil.
In each case, I show that they are the only $P$-equivariant rational
map $T^*_x G/P \dasharrow (G/Q)_x$
(propositions \ref{unicite_spin} and \ref{unicite_sl5}),
if $(G/Q)_x$ denotes the variety of $y$'s in $G/Q$ with stabilizor
$Q_y$ such that $P_x \cap Q_y$ is parabolic.

In the case of a flop
of type $E_{6,I}$ for instance, we get a $Spin_{10}$-equivariant
rational map; $G/P$ is often called the Moufang plane $\p^2_\O$
(it is some kind of octonionic projective plane). As an example of the
above discussion, the restriction of the flop to a cotangent space
should interpret as the
``quotient map'' $\A^2_\O \dasharrow \p^1_\O$. In the first
section, I show that this map has some properties of such a quotient;
for example, its
fibers carry a natural structure of algebra
isomorphic with the octonions (corollary \ref{fibre_octave}).
I also study the projective geometry of the corresponding spinor
variety.

Similarly, section \ref{sl5} gives a model for the restriction of the
Mukai flop of type $E_{6,II}$ to a cotangent space. In this case, a
Levi factor contains $SL_2 \times SL_5$ and the relevant factor of
$T^* E_6/P_3$ is $Hom(\C^2,{(\wedge^2 \C^5)}^*)$. The given classification
of $(GL_2 \times GL_5)$-orbits in
$Hom(\C^2,{(\wedge^2 \C^5)}^*)$ allows to 
understand the $E_6$-orbits in $T^* E_6/P_3$. Finally, corollary
\ref{lieu} states that the Mukai flop is defined only on the open
orbit of $T^*E_6/P_3$, and describes the image of all orbits in
$T^*E_6/P_3$ as
nilpotent orbits in $\mathfrak{e}_6$.

\lpara
{\bf Acknowledgement } I thank Baohua Fu for many usefull discussions on
the topic of nilpotent orbits, and stimulating questions.

\tableofcontents

\sectionplus{The group $Spin_{10}$ and $\p^1_\O$.}

\label{spin10}

\subsectionplus{Geometric definition of composition algebras}

For more details on composition algebras, the reader my consult 
\cite{projective}. I
recall that if $R$ is a ring, then $\RR,\CR,\HR$ and $\OR$ denote the
four usual split composition algebras over $R$. Therefore, $\RR = R,
\CR = R \oplus R$, $\HR$ is the algebra of $2 \times 2$-matrices with
coefficients in $R$, and $\OR$ is obtained from $\HR$ by
Cayley-Dickson's process.

Their norms will be
denoted $N$. If $\a$ is one of those and 
$z \in \a$, then
$L_z$ and $R_z$ denote the endomorphisms of $\a$ of left and right
multiplication by $z$, and $L(z),R(z) \subset \a$ their images.

In the following, we will have to define a composition algebra
structure on a vector space by geometric means. This subsection
explains how it is possible.
In this section, $k$ is an algebraically closed field of 
caracteristic different from 2.

\begin{prop}
\label{composition_associative}
Let $V$ be a
$k$-vector space of dimension $a$, with $a \in \{1,2,4\}$. 
Let $x_0 \in V - \{0\}$ 
and ${\cal N} \subset \p V$ be a smooth quadric such
that the class of $x_0$ in $\p V$ does not belong to $\cal N$.

Then if $a \in \{1,2\}$, there exists a unique composition algebra
structure on $V$ with unit $x_0$ and such that ${\cal N}$ is the
quadric of elements with vanishing norm. If $a=4$, there are two such
composition algebras.
\end{prop}
\noindent
Therefore, giving a composition algebra structure on $V$ is
equivalent with giving a smooth quadric in $\p V$ and a point out of
its affine cone (if $a=4$, we must moreover choose a component of the
variety of maximal isotropic subspaces of the quadric).\\
\pr
The existence of the algebra is an immediate consequence of the fact
that $Aut({\cal N})$ acts transitively on $\p V - {\cal N}$.

The unicity in the cases when $a \in \{1,2\}$ is easy. Assume $a=4$ and
let $(x,y) \mapsto xy$ be a composition product satisfying the
conditions of the proposition. Let 
${\cal L},{\cal R} \simeq \p^1$ be
the two families of isotropic lines. For $x \in {\cal Q}$, denote
$l(x)$ (resp. $r(x)$) the isotropic line in $\cal L$ (resp. $\cal R$)
containinig $x$. Up to changing the algebra
structure $(x,y) \mapsto xy$ into $(x,y) \mapsto yx$, we may assume
that $\forall x \in {\cal N},L(x) \in {\cal L}$. Therefore,
$L(x)=l(x)$ and $R(x)=r(x)$.

If $z \in V$, let $[z]$ denote its class in $\p V$.
Then, for generic $x,y \in \cal N$, $[xy] = L(x) \cap R(y) = l(x) \cap r(y)$.
Therefore, the product of two elements of $\cal N$ is fixed up to a
scale once 
$\cal N$ is. One checks also that $xy=0$ \iff $r(x) \cap l(y)$ is
orthogonal to the unit, with respect to the scalar product defined by
$\cal N$. In view of lemma \ref{app_lin} applied to left
multiplication by $x \in {\cal N}$, the proposition is proved.
\qed

\begin{lemm}
Let $V$ and $W$ be vector spaces and $f,g : V \rightarrow W$ linear maps. 
Let $X \subset \p V$ be an irreducible variety included in
no hyperplane of $\p V$. Assume that the induced rational maps 
$[f]_{|X},[g]_{|X} : X \dasharrow \p W$ are equal and that 
$\ker f = \ker g$. Then there exists $\lambda \in k-\{0\} : f = \lambda g$.
\label{app_lin}
\end{lemm}
\pr
$f$ and $g$ have the same image, spanned by $f(X)=g(X)$. They also
have the same
kernel by hypothesis. Therefore, there is a linear automorphism $h$ of
this common image, such that $g=hf$. Since $[f]_{|X} = [g]_{|X}$, any
vector in $f(X)$ is an eigenvector for $h$, from which the lemma
follows.
\qed

\lpara

We now consider the case of the octonions.

\begin{prop}
Let $V$ be an 8-dimensional vector space and ${\cal N} \subset \p V$ a
smooth quadric. Let $G$ denote the grassmannian of maximal isotropic
subspaces of $\cal N$, and let $l$ be an
isomorphism between $\cal N$ and an irreducible component of $G$.
Assume $\forall x \in {\cal N},x \in l(x)$ and let
$x_0 \in V - \{ 0 \}$ such
that $[x_0] \not \in {\cal N}$.

Then there exists a unique composition algebra structure on $V$ with
unit $x_0$ and such
that for all $x \in {\cal N}$, we have $l(x)=L(x)$.
\label{octave}
\end{prop}
\pr
Given an octonionic structure on $V$,
it is known as ``triality principle'' \cite[chapter IV]{chevalley}
that
$L$ is an isomorphism on its image,
which is a connected component of $G$.

The unicity of the algebra structure follows the lines of the previous
proposition. Let $(x,y) \mapsto xy$ be an algebra structure on $V$
with unit $x_0$ and such that $L=l$. If $x \in {\cal N}$ is generic,
then the line
$(x,x_0)$ meets ${\cal N}$ at $x$ and $\overline x$. Therefore, $x_0$
determines the conjugation. By hypothesis, $L(x)=l(x)$, therefore
we get $R(x)$ as
$\overline{l(\overline x)}$. Now, the class of the product $xy$ 
in $\p V$ is again 
$L(x) \cap R(y)$, and $xy=0$ \iff $\dim L(x) \cap R(x) = 3$, as is
well-known \cite[IV.4.2]{chevalley}. Therefore, lemma \ref{app_lin}
proves the unicity of the algebra.

Let us prove its existence. Put on $V$ an arbitrary
structure of composition algebra $(x,y) \mapsto xy$
such that the quadric of elements of vanishing norm is $\cal N$. 
This induces
isomorphisms $L,R$ between ${\cal N}$ and the components of $G$. Set
$r(x) = \overline{l(\overline x)}$. We can assume that $L$ and $l$
have the same image. Therefore, there exist $f,g \in Aut({\cal N})$
such that $l(x)=L(f(x))$ and $r(x)=R(g(x))$.
The hypothesis $x \in l(x)$ implies 
$x \in r(x)$, and so $f(x) \in L(x)$ and $g(x) \in R(x)$
\cite[proposition 1.1]{projective}. By the
following lemma \ref{fois_octave}, there exist invertible $\alpha,\beta$ such
that $f(x)=x\alpha$ and $g(y)=\beta y$. The composition algebra
$x*y=(x\alpha)(\beta y)$, with unit $\beta^{-1} \alpha^{-1}$, 
satisfies the conditions of the proposition.
\qed

\begin{lemm}
Let $m:\ok \rightarrow \ok$ a linear map preserving $\cal N$ and such that 
$\forall x \in {\cal N},m(x) \in L(x)$. Then there exists $\alpha \in \ok$
such that $\forall x \in \ok,m(x)=x\alpha$.
\label{fois_octave}
\vspace{-.4cm}
\end{lemm}
\pr
Left to the reader \cite[p.48]{these}.
\qed

\subsectionplus{The 8-dimensional quadric as $\p^1_\O$}

\label{p1o}

I have juste recalled the triality principle, which implies that
the three 8-dimensional
fundamental representations of $Spin_8$ can be identified with the
algebra of octonions.
The goal of this subsection is to relate the group $Spin_{10}$ with the
octonions, see proposition \ref{s+-}. To study the
representations of $Spin_{10}$, my strategy is to restrict them to
representations of $Spin_8$.
Before proving proposition \ref{s+-},
I need to make a computation in Clifford algebras. My notations are 
those of \cite{chevalley}.
\lpara

Let $V$ be a $k$-vector space of even dimension and 
equipped with a non-degenerate
quadratic form $q$. Let $V' \subset V$ be a codimension
two subspace in $V$ such that $q_{|V'}$ is non-degenerate. 
Let $C,C'$ denote the Clifford algebra of $V,V'$
(the Clifford algebra of $V$ is the tensor algebra of $V$ mod out by
the relations $x \otimes x = q(x)$). Let $\alpha$ be the ``main
antiautomophism'' of $V$, defined by 
$\alpha(v_1 \ldots v_k)=v_k \ldots v_1$.

Let 
$V' = N' \oplus P'$ be a decomposition into isotropic subspaces. Let
$x_0,y_0 \in V$ be orthogonal to $V'$ and such that
$q(x_0,y_0)=1$. Denote $N = N' \oplus k.x_0$ and $P = P' \oplus k.y_0$.

Let $C_N \subset C$ (resp. $C_N' \subset C'$) be the
subalgebra of $C$ (resp. $C'$) generated by $N$ (resp. $N'$). Let 
$f' \in C_N'$ be the product of the elements of a basis of $N'$ and
$f=f' y_0$. Let $\S^\pm$ and ${\S'}^\pm$ be the spin representations of
$Spin(V)$ and $Spin(V')$. We may choose $\S^+$ (resp. $\S^-$) be the
subspace of even (resp. odd) elements of $C_N$, and similarly for
${\S'}^\pm$.
 
There are isomorphisms 
$\varphi^\pm$ between ${\S'}^+ \oplus {\S'}^- = C_N'$
and $\S^\pm = C_N^\pm$, given by $\varphi^+(u'_++u'_-) = u'_+ + u'_-x_0$
and $\varphi^-(u'_++u'_-) = u'_+x_0 + u'_-$.
Finally, there is a quadratic
map $\beta : C_N \times C_N \rightarrow \wedge V$,
where $\beta(u,v)$ is the image of $uf\alpha(v) \in C$ in $\wedge V$
under the
canonical vector space isomorphism $C \simeq \wedge V$ 
\cite[p.102,103 and II 1.6]{chevalley}. Let 
$\beta' : C_N' \times C_N' \rightarrow \wedge V'$ be the similar map for $V'$.

\begin{prop}
Let $r' = \dim V'/2$.
Let $u'_+,v'_+ \in {\S'}^+$ and let $u'_-,v'_- \in {\S'}^-$. We have
$$
\begin{array}{cl}
& \beta[\varphi^+(u'_++u'_-),\varphi^+(v'_++v'_-)] \\
= & \beta'(u'_+,v'_+) \wedge y_0
- x_0 \wedge y_0 \wedge \beta'(u'_+,v'_-) + \beta'(u'_+,v'_-) \\
+ & (-1)^{r'}(x_0 \wedge y_0 \wedge \beta'(u'_-,v'_+) +  \beta'(u'_-,v'_+))
- x_0 \wedge \beta'(u'_-,v'_-),
\end{array}
$$
and
$$
\begin{array}{cl}
& \beta[\varphi^-(u'_++u'_-),\varphi^-(v'_++v'_-)] \\
= & x_0 \wedge \beta'(u'_+,v'_+)
+ (-1)^{r'}(x_0 \wedge y_0 \wedge \beta'(u'_+,v'_-) + \beta'(u'_+,v'_-)) \\
+ & y_0 \wedge x_0 \wedge \beta'(u'_-,v'_+) +  \beta'(u'_-,v'_+)
- \beta'(u'_-,v'_-) \wedge y_0.
\end{array}
$$
\label{beta}
\end{prop}
\pr
We have, in the Clifford algebra $C$, 
$u'_+ f' y_0 \alpha(v'_+) = u'_+ f' \alpha(v'_+)y_0$, so
$\beta[\varphi^+(u'_+),\varphi^+(v'_+)] = \beta'(u'_+,v'_+) \wedge y_0$.
We can compute the other terms
$\beta[\varphi^+(u'_\pm),\varphi^+(v'_\pm)]$ using the facts
$$
\begin{array}{l}
u'_+f'y_0\alpha(v'_-x_0)=y_0x_0u'_+f'\alpha(v'_-), \hspace{.25cm}
u'_-x_0f'y_0\alpha(v'_+) = (-1)^{r'} x_0y_0u'_-f'\alpha(v'_+), \\
\mbox{and }\ u'_-x_0f'y_0\alpha(v'_-)=x_0y_0x_0u'_-f'\alpha(v'_-)
=x_0u'_-f'\alpha(v'_-).
\end{array}
$$

The computation of $\beta[\varphi^-(u'_++u'_-),\varphi^-(v'_++v'_-)]$
is similar.
\qed

\lpara

Our second task is to describe spinor representations using
octonions. Let $V=H_2(\ok)$ denote the 10-dimensional
$k$-vector space of $2 \times 2$
hermitian matrices with entries in $\ok$. Let $\det$ be the quadratic
form on $H_2(\ok)$ defined by
$\det\left (\matdd tz{\overline z}u \right ) = tu - N(z)$ ($t,u \in k$
and $z \in \ok$). Recall \cite[III 1.2,III 1.4]{chevalley}
that the variety of maximal
isotropic subspaces of $V$ has two components; they will be denoted
$G_Q^+(5,V)$ and $G_Q^-(5,V)$. Moreover, there are natural projective
embeddings $G_Q^\pm(5,V) \subset \p \S^\pm$ in the projectivized
spinor representations, the elements of $\S^\pm$ which class are in 
$G_Q^\pm(5,V)$ being called ``pure spinors''.

Let $\nu_2^+:\ok \times \ok \rightarrow H_2(\ok)$ the
quadratic map defined by $\nu_2^+(a,b)=\matddr{N(a)}{a\overline b}
{b \overline a}{N(b)}$ and $\mu^+$ the polarization of $\nu^+\ :\ 
\mu^+((a,b),(c,d))=\nu_2^+(a+c,b+d)-\nu_2^+(a,b)-\nu_2(c,d)$).
Similarly, let $\nu_2^-(a,b)=\matddr{N(b)}{a\overline b}
{b \overline a}{N(a)}$ and $\mu^-$ the polarization of $\nu_2^-$.

Let $X^+=X^- \subset \p(\ok \oplus \ok)$ be defined by
$[(a,b)] \in X^\pm \Longleftrightarrow \nu_2^\pm(a,b)=0$.

\begin{prop}
The variety $X^\pm$ is isomorphic with $G_Q^\pm(5,V)$.
An isomorphism $X^\pm \rightarrow G_Q^\pm(5,V)$ maps $(u,v)$ on
the image of $\mu^\pm((u,v),.)$.
\label{s+-}
\end{prop}
\pr
Let $q=-\det$,
$V' = \left \{ \matddr 0z{\overline z}0 \right \} \simeq \ok$,
$x_0=\matddr 000{-1}$ and $y_0=\matddr 1000$.
Let $\beta_k$ denote the component in $\wedge^k V \subset \wedge V$ of
$\beta$.
Since $q$ restricts to the norm of octonions on $V' \simeq \ok$,
by the triality principle \cite[Chapter IV]{chevalley}, 
with the notations of proposition \ref{beta}, there are
linear isomorphisms ${\S'}^\pm \rightarrow \ok$ such that the
map $\beta'_1 : {\S'}^+ \times {\S'}^- \rightarrow V'$ identifies with
the product of octonions, and 
${(\beta'_0)}^+ : \S^+ \times \S^+ \rightarrow k, 
{(\beta'_0)}^- : \S^- \times \S^- \rightarrow k$ identify with the
scalar product of octonions.
Composing with the automorphism 
$b \mapsto \overline b$, of ${\S'}^- \simeq \ok$, we may assume that
$\beta'_1$ is in fact given by $(a,b) \mapsto a \overline b$.

By proposition \ref{beta}, $\S^+$ and $\S^-$ therefore identify with
$\ok \oplus \ok$ in such a way that
$\beta_1^+((a,b),(a,b))=2N(a)y_0-2N(b)x_0+2a\overline b$ and
$\beta_1^-((a,b),(a,b))=-2N(a)x_0+2N(b)y_0+2a\overline b$, that is to say,
$\beta_1^\pm = \mu^\pm$.

By proposition \cite[III 5.2]{chevalley} the spinor varieties 
$G_Q^\pm(5,V) \subset \p \S^\pm$ are defined by the
equations $N(a)=N(b)=0,a\overline b=0$, which is equivalent to
$\nu_2^\pm = 0$. Therefore, they are isomorphic with
$X^\pm$. Moreover, since the linear space corresponding to $s$ is the
image of $\mu^\pm(s,.)$ \cite[III 4.4]{chevalley}, the proposition is proved.
\qed

\lpara

In the sequel, we will identify both $\S^+$ and $\S^-$ with $\ok \oplus \ok$,
keeping however in mind the fact that $\S^+$ and $\S^-$ are non-equivalent
$Spin_{10}$-representations.
The projectivization $\p \S^\pm$ of $\S^\pm$ have two
$Spin_{10}$-orbits, by \cite[prop. 2 p.1011]{igusa}. 
The closed orbits are $X^+$ and $X^-$.

\lpara

Now comes the explanation of the title of this subsection : the variety
of classes of matrices 
$\left [ \matdd tz{\overline z}u \right ]\in \p V$ with $tu-N(z)=0$ 
is a $Spin_8$-conformal
compactification of the variety of classes of matrices of the form
$\left [ \matdd 1z{\overline z}{N(z)} \right ]$ which is isomorphic with 
$\ok \simeq \A^1_\O$, therefore, it can be thought as
$\p^1_\O$. Moreover, the projectivisations
$\overline{\nu_2^\pm} : \S^\pm \dasharrow \p\{\det = 0\}$ of
the maps $\nu_2^\pm$ are some kind of quotient
maps $\A^2_\O \dasharrow \p^1_\O$. Proposition \ref{qsinterqt} and
corollary \ref{fibre_octave} illustrate this viewpoint.

For the moment, we show that $\overline {\nu_2^+}$ and
$\overline {\nu_2^-}$ are the only natural (ie $Spin_{10}$-equivariant)
candidates for such a kind of quotient (proposition \ref{unicite_spin}).
Let $Q \subset \p V$ denote the quadric defined by $\det$.

\begin{lemm}
There is a unique 15-dimensional $Spin_{10}$-orbit in
$(\p \S^+ - X^+) \times Q$.
\vspace{-.4cm}
\label{unique_dim15}
\end{lemm}
\pr
Let $(s_1,x_1),(s_2,x_2) \in (\p \S^+ - X^+) \times Q$. We may assume
that $s_1=s_2=s$. Let $G_0 \subset Spin_{10}$ be the stabilizor of $s$.
From the proof of \cite[prop. 2 p.1011]{igusa}, it follows that
$\overline{\nu_2^+}(s) \in Q$ is the only line in $Q$ stabilized by $G_0$.
Therefore, $x_1 = x_2 = \nu_2^+(s)$.
\qed

\begin{prop}
$\overline {\nu_2^+} : \S^+ \dasharrow Q$ is the only 
$(k^* \times Spin_{10})$-equivariant rational map
$\S^+ \dasharrow Q$.
\label{unicite_spin}
\end{prop}
\pr
Let $\nu : \S^+ \dasharrow Q$ be any $(k^* \times
Spin_{10})$-equivariant rational map
$\S^+ \dasharrow Q$. Then $\overline {\nu_2^+}$ and $\nu$ induce rational
maps $\p \S^+ \dasharrow Q$, which will be denoted with the same
letter. Since $\nu$ is $Spin_{10}$-equivariant, it is
defined on $\p \S^+ - X^+$. Therefore, the variety of
$\{(s,\nu(s)) : s \in \p \S^+ - X^+\}$ is a 15-dimensional orbit in 
$(\p \S^+ - X^+) \times Q$. By
lemma \ref{unique_dim15}, it is equal to the orbit
$\{(s,\overline {\nu_2^+}(s)) : s \in \p \S^+ - X^+\}$.
\qed

\subsectionplus{Projective geometry of the spinor variety}

We keep the notations of the previous subsection; namely,
$V=H_2(\ok)$, $\S^+ = \S^- = \ok \oplus \ok$ are the two spinor
representations of $Spin_{10}$, and $\nu_2^\pm : \S^\pm \rightarrow V$ are
the quadratic $Spin_{10}$-equivariant maps defined above. Their
polarizations are denoted $\mu^\pm$. We
denote $Q \subset \p V$ the smooth quadric defined by $\det$.
If $(a,b) \in \ok \oplus \ok$ we denote $[a,b]$ its class in
$\p (\ok \oplus \ok)$. Finally, if $X \subset \p^n$ is a variety and
$x \in X$, let $T_xX$ its tangent space
and let $\widehat X \subset \A^{n+1}_k$ denote the affine cone over $X$.

\lpara

Recall from \cite[III 2.3]{chevalley} that there is a 
$Spin_{10}$-equivariant perfect pairing
$\S^+ \times \S^- \rightarrow k$.
This allows identifying $\S^-$ with the dual of $\S^+$. Recall that the
dual variety of a variety $X$ is the closure of the set of tangent hyperplanes,
where a tangent hyperplane is by definition a hyperplane containing a
tangent space $T_xX$ at a smooth point $x \in X$.

\begin{prop}
The equivariant isomorphism $\p {\S^+}^* \simeq \p \S^-$ identifies
the dual variety of $X^+$ with $X^-$.
\end{prop}
\pr
The dual variety of $X^+$ is a $Spin_{10}$-stable closed variety. Since in
$\p S^-$ there are only two orbits, by \cite[prop 2 p.1011]{igusa}, 
it is either the
whole projective space $\p \S^-$, which is absurd, or the variety
$X^-$.
\qed

\lpara

If $X \subset \p^n$ is a subvariety of projective space, and if
$z \in \p^n - X$, the entry locus of $z$ is classically defined
as the
closure of the set of points $x \in X$ such that the line joining $x$
and $z$ meets $X$ at at least two distinct points.

If $s \in \p \S^\pm - X^\pm$, denote $L_s^\pm$ the variety 
$\overline{{(\nu_2^\pm)}^{-1}(k^*.\nu_2^\pm(t))} \subset \S^\pm$,
where $t \in \S^+$ is
such that $[t]=s$. Let $L^\pm$ denote the
variety $\{(s,v) \in (\p \S^\pm - X^\pm) \times \S^\pm : v \in L_s^\pm \}$.
Finally, let 
$\overline \nu_2^\pm : \p \S^\pm - X^\pm \rightarrow Q$ denote the map
induced by $\nu_2^\pm : \S^\pm \rightarrow V$.

\begin{prop}
Let $s \in \p \S^+ - X^+$. Then the entry locus $Q^+_s$ of $s$ in 
$X^+$ is a smooth 6-dimensional
quadric in the 7-dimensional projective space $\p L_s^+$. Moreover,
the fibration $L \rightarrow \p \S^+ - X^+$ is locally trivial and is the
push-back by $\overline \nu_2^+$ of a vector bundle on $Q \subset \p V$.
\label{ls}
\end{prop}
\rek The bundle over $Q$ of the proposition is often called the
spinor bundle.\\
\pr
Since $\p S^+ - X^+$ is a single $Spin_{10}$-orbit, it is enough to
check the first claim of the
proposition for $s=[1,0]$. Computing $Q^+_s$ is equivalent
with solving the equation $(1,0)=(a,x)+(b,y)$ in $\ok \oplus \ok$,
with $(a,x)$ and $(b,y)$ in the affine cone over $X^+$. Equivalently,
$(a,x)$ satisfies $N(a)=N(x)=0$ and $a \overline x=0$, and similarly
for $(b,y)$.

Now, the equality $a+b = 1$ implies $N(a,b) \not = 0$ ($N(.,.)$
denotes the polarization of $N$). This, in turn, implies 
$R(a) \cap R(b) = \{ 0 \}$ \cite[IV 4.4]{chevalley}. 
Since $a \overline x = 0$,
$x \in R(a)$ \cite[proposition 1.1]{projective}
and similarly $y \in R(b)$. Since $x = -y$,
it follows that $x \in R(a) \cap R(b)$, so $x=0$.

We thus have proved that the entry locus $Q^+_s$ is included in the
variety of elements $[a,0]$ with $N(a) = 0$. Conversely, this smooth
quadric is included in $Q^+_s$.
Since if $N(a) \not = 0$, then left multiplication by $a$ is
invertible, and
a direct computation shows that $L_s^+= \ok \oplus \{0\}$.

\lpara

To show that $L^+$ is a vector bundle, let
$s \in \S^+$ and $x=[\nu_2^+(s)] \in Q$. Let $L_x \subset V$ denote
the line corresponding to $x$.
First recall by definition 
that the image of the restriction of $\nu_2^+$ to $L_s^+$ is
the line of multiples of $\nu_2^+(s)$. Therefore, 
$\mu^+(s,L_s^+)=k.\nu_2^+(s)$.
The linear space $\mu^+(s,\S^+)$ is $\widehat{T_xQ}$,
thus the kernel of the composition
$\S^+ \stackrel{\mu^+(s,.)}{\longrightarrow} \widehat{T_xX} 
\rightarrow \widehat{T_xX}/L_x$ is exactly
$L_s^+$. Therefore, $L^+$ is the kernel
of a morphism of vector bundles over $\p \S^+ - X^+$ with constant
rank; so, it is locally free.

Since $L_s^+$ is constant on a fiber 
${(\overline \nu_2^+)}^{-1}(x)$, $L^+$ is the
push-back of a vector bundle on $Q$.
\qed

\lpara

We now study the family of quadrics $\{Q^+_s\}$. Let $G(8,\S^+)$ denote
the grassmannian of 8-dimensional linear spaces in $\S^+$ and consider
the variety ${\cal Q} \subset G(8,\S^+)$ of 8-dimensional linear
spaces $L$ in $\S^+$ such that $X^+ \cap \p L$ 
is a smooth 6-dimensional quadric.
\begin{prop}
The variety $\cal Q$ is $Spin_{10}$-equivariantly
isomorphic with $Q$ and any element of $\cal Q$ is of the form $Q_s$.
Moreover, let 
$s,t \in \p \S^+ - X^+$; one of the following holds~:
\begin{enumerate}
\item
$L_s^+=L_t^+$ and $Q_s=Q_t$.
\item
$Q_s^+ \cap Q_t^+ = \p L_s^+ \cap \p L_t^+ \simeq \p^3$.
\item
$\p L_s^+ \cap \p L_t^+ = \emptyset$.
\end{enumerate}
\label{qsinterqt}
\end{prop}
\rek
Although this does not make sense due to the lack of associativity of
the octonions, the maps $\nu_2^\pm : \S^\pm \dasharrow Q$ 
should be some kind of
quotient maps $\A^2_\O \dasharrow \p^1_\O$. The linear space $L_s^+$ can
be interpreted as the set of $\ok$-multiples of $s$ (in $\ok \oplus \ok$). 
With this point of view, the proposition says that for two
non-degenerate (out of $X^+$) vectors in 
$\ok \oplus \ok$, there are three possibilities : either they are
linked (1), either they are free (3), either they are ``weakly linked''
(2). This last case would not occur if we would consider non-split
octonions, say for example over the field $\R$ of real numbers.
The same situation holds when one considers two
non-degenerate vectors $v,w \in \hk \oplus \hk$. 
In fact, it is easy to check that
$$\dim\ (\{ \lambda.v : \lambda \in \hk \} \cap
\{ \lambda.w : \lambda \in \hk \} ) \in \{ 0,2,4 \}.$$ 
(see the remark after lemma 2.1
in \cite{projective}). \\
\pr
Proposition \ref{ls} and the fonctorial property of
grassmannians show that there is a map 
$\psi : Q \rightarrow {\cal Q}$.
In the other way, let $l\in \cal Q$. 
Let $\delta$ be a generic line in $\p l$; this line meets the quadric
$X^+ \cap \p l$ in two points $x$ and $y$. Since $\nu_2^+$ vanishes on 
$\widehat X^+$, for any $s$ in
$\delta$, we have 
${\overline \nu}_2^+(s)={\overline \mu}^+(x,y)$. 
Therefore, ${\overline \nu}_2^+$ is constant on the generic
lines in $\p l$, so it is constant on $\p l$.
This proves that there is a map 
$\varphi : {\cal Q} \rightarrow Q$, induced by $\nu_2^+$. 

It is obvious, by construction, that $\varphi$ and $\psi$ are inverse
maps, so the first point of the proposition is proved.

The rest of the proposition follows. In fact,
set $s=(1,0)$, so that $L_s^+ = \ok \oplus 0$.
If $t=s$, then $L_s^+ = L_t^+$.
If $t=(1,b)$, with $N(b)=0$, then an easy computation shows that
$\p L_s^+ \cap \p L_t^+ = Q_s^+ \cap Q_t^+ = \{ (c,0) : c \in R(b) \}$.
If $t=(0,1)$, then $L_t^+ = 0 \oplus \ok$ and so
$\p L_s^+ \cap \p L_t^+ = \emptyset$.

Since there are three 
$Spin_{10}$-orbits in $Q \times Q$, these three
examples exhaust all the possibilities for a couple 
$(L_s^+,L_t^+) \in {\cal Q} \times {\cal Q}$.
\qed
\lpara

Let $s \in \p \S^+ - X^+$. Define $Q^-_s$ as the intersection of $X^-$
with the orthogonal of $L_s^+$ (in other words, $Q^-_s \subset {(X^+)}^*$ 
is the variety of tangent hyperplanes which contain $L_s^+$).

\begin{prop}
With notations above, $Q_s^-$ is a 6-dimensional smooth quadric in
$X^-$. Moreover, its linear span in $\S^-$ is the closure of
${(\nu_2^-)}^{-1}(k^*.\nu_2^+(s))=:L_s^-$.
\vspace{-.4cm}
\label{qs-}
\end{prop}
\pr
Arguing as in the proof of proposition \ref{s+-} one can show
that the equivariant duality
between $\S^+$ and $\S^-$ is 
$\scal{(a,b),(c,d)} = N(a,c) + N(b,d)$.
Therefore, if $s=[1,0]$, then $Q_s^-$ is the variety
$\{[0,b] : N(b)=0\}$. Its linear span is $0 \oplus \ok$, which is sent
by $\nu_2^-$ on $\matddr 1000 = \nu_2^+(s)$.
\qed

Let $\varphi^\pm$ denote the isomorphisms between $X^\pm$ and the
components of the grassmannian of maximal isotropic subspaces in $Q$.
We have another caracterisation of the quadrics $Q_s^+$ and $Q_s^-$~:

\begin{prop}
Let $x \in X^\pm$ and $s \in \S^+$. Then $x \in Q_s^\pm$ \iff
$\nu_2^+(s) \in \varphi^\pm(x)$.
\label{qs}
\end{prop}
\pr
By proposition \ref{ls}, there exists a quadratic form 
$q_s$ on $L_s^+$, which zero locus is
$Q_s^+$, and such that 
$\forall u \in L_s^+, \nu_2^+(u)=q_s(u) \nu_2^+(s)$.
Let $x \in \widehat Q_s^+ \subset L_s^+$, then for $u \in L_s^+$, we have
$\mu^+(x,u) = q_s(x,u) \nu_2^+(s)$.
Therefore,
$\nu_2^+(s) \in \{\mu^+(x,u):u\in \S^+\} = \varphi^+(x)$.

The converse implication 
$\nu_2^+(s) \in \varphi^+(x) \Rightarrow x \in Q_s^+$
follows by a dimension count argument.

In view of proposition \ref{qs-}, the proof of the same result for
$Q_s^-$ is similar.
\qed

\subsectionplus{Equivariant octonionic structure on the fibers on $\nu_2^\pm$}

In a honest projective space $\p^n_k$, over a field $k$, the choice of
an element $v \in \A^{n+1}_k$ identifies the closure of
the fiber of the quotient map 
$\A^{n+1}_k \dasharrow \p^n_k$ with $k$, since any element in this
fiber can uniquely be written as $\lambda.v$, with $\lambda \in k$. 
Therefore, this fiber
carries the structure of a field, isomorphic with $k$.

We will see (corollary \ref{fibre_octave})
in this subsection something analogous for $\nu_2^\pm$,
which is interpreted as a quotient map. However, let $s \in \S^+$;
the image of the stabilizer of $s$ in $GL(L_s^+)$ contains $Spin_7$ by
\cite[prop 2 p.1011]{igusa}. 
Therefore, there is no hope to give $L_s^+$ an equivariant
octonionic structure.

I will show that given two generic spinors $s,t \in \S^+$, there are
equivariant octonionic structures on $L_s^+$ and $L_t^+$ (and indeed the
stabilizor of two elements has a quotient isomorphic with $G_2$). 
I don't know how
to interpret the necessity of two spinors to define such a structure
in terms of octonionic projective geometry.

\lpara

Let $s,t \in \p \S^+ - X^+$ such that 
$\scal{\nu^+_2(s),\nu^+_2(t)} \not = 0$.
The idea of the geometric definition of an octonionic structure on $L_s^+$
is as follows~: we have the two quadrics $Q_s^+$ and
$Q_s^-$. Let $Q_s$ denote variety of lines in $Q$
containing $[\nu_2^+(s)]$. Then $Q_s$ is isomorphic with a
6-dimensional quadric.
By proposition \ref{qs}, 
$Q_s^+$ and $Q_s^-$ parametrize the maximal
isotropic linear spaces of $Q_s$. The point is to
show that $s$ and $t$ yield an isomorphism
$Q_s^+ \stackrel{\sim}{\rightarrow} Q_s^-$. 
Then, proposition \ref{octave} gives the
octonionic structure.

The next proposition yields the isomorphism 
$\psi:Q_s^+ \stackrel{\sim}{\rightarrow} Q_s^-$.
Let $x \in Q_s^+$ be such that the line through $x$ and $s$ is not a
tangent line to $Q_s^+$. Call $\overline x$ the other point of intersection of
this line with $Q_s$. Moreover, set 
$r(x)={\scal{T_xX^+,L^+_s}}^\bot \subset \p {\S^+}^* = \p \S^-$.

\begin{prop}
For all $x \in Q_s^+$, $r(x)$ is a maximal isotropic subspace of $Q_s^-$.
Moreover, if $(x,s)$ is not a tangent line to $Q_s$, then $r(x)$ and 
$r(\overline x)$ are supplementary subspaces of $L_s^-$. 
Call $\psi(x)$ the image of $t$
by the projection on $r(x)$ with center $r(\overline x)$. Then 
$\psi : Q_s^+ \rightarrow Q_s^-$ is an isomorphism.
\label{iso}
\end{prop}
\pr
Assume $s=[1,0]$ and $t=[0,1]$. Let $x=[a,0] \in Q_s^+$ 
(therefore $N(a)=0$). Since $X^+$ is the variety of pairs $[a,b]$ with
$N(a)=N(b)=0$ and $a \overline b = 0$, 
$T_xX^+=\{[c,d]:N(a,c)=0 \mbox{ and } a \overline d = 0\}$. Therefore,
its orthogonal is the set of $[c,d]$ with $c$ colinear with $a$ and
$d \in R(c)$. So $r(x)=\{[0,d]:d \in R(a)\}$. This is endeed a maximal
isotropic subspace of $Q_s^-$.

Moreover, we have $\overline x = [\overline a,0]$, and so
$r(\overline x) = \{ [0,d] :d \in R(\overline a) \}$. Therefore, $r(x)$
and $r(\overline x)$ are supplementary.

Finally, since $t=[0,1]=[0,(a+\overline a)/2]$, we deduce that
$\psi(x)=[0,a]$. We have proved that
$\psi([a,0])=[0,a]$, so $\psi$ is an isomorphism.
\qed

\begin{coro}
Let $s,t \in \p \S^+ - X^+$ such that 
$\scal{\nu^+_2(s),\nu^+_2(t)} \not = 0$.
Then $L_s^+$ has a natural structure of algebra, isomorphic with $\ok$.\\
Moreover, when $(s,t)$ vary, this octonionic structure on the vector
bundle with fiber $L_s^+$ varies algebraically.
\label{fibre_octave}
\end{coro}
\pr
We have isomorphisms $\psi^\pm$ between $Q_s^\pm$ and the components of
the variety of maximal isotropic subspaces in $Q_s$, as explained at
the beginning of this paragraph.

If $s=(1,0)$ and $t=(0,1)$, it follows from the proof of
proposition \ref{iso} that 
$\forall x \in Q_s^+, \dim\ ( \psi^+(x) \cap \psi^-(\psi(x)) ) = 3.$
This is analogous to the condition $x \in l(x)$
of proposition \ref{octave}, and
therefore $s$ and the isomorphisms $\psi,\psi^+,\psi^-$ define a
unique octonionic structure on $L_s^+$.

\lpara

It follows by general arguments that this octonionic structure varies
algebraically. Alternatively, one can give another construction of
this octonionic structure, where the algebraicity is clear.

Let $x=\nu_2^+(s)$ and $M = \widehat {T_xQ}/k.x$.
Then, as one checks one the example $s=(1,0),t=(0,1)$,
$\mu^+(s,.)$ restricts to an isomorphism $\nu_t$ between $L_t^+$ and $M$ and
$\mu^+(t,.)$ to an isomorphism $\nu_s$ between $L_s^+$ and $M$. 
We can therefore give an octonionic structure to $L_s^+$ by setting
$$
\forall u,v \in L_s^+, \ uv = \nu_s^{-1}[\mu^+(u,\nu_t^{-1}(\nu_s(v)))].
$$
A direct computation shows that this octonionic structure is the same
as the previous one.
\qed

%******************************************************************************

%******************************************************************************

%******************************************************************************

\sectionplus{Geometry associated with two skew-forms in $k^5$}

\label{sl5}

In this section, we consider a model for the restriction of the Mukai
flop of the second kind to a tangent space. In the first subsection, 
we prove
lemmas which will suffice defining this restriction, in section
\ref{flop}. The second subsection will be used when classifying the
orbits in $T^*G/P$, for $G$ of type $E_6$ and $P$ corresponding to
$\alpha_3$. The third subsection shows that the involved rational map
is the unique equivariant rational map.

Let $k$ denote an arbitrary field.

\subsectionplus{A rational map
$Hom(k^2,{(\wedge^2 k^5)}^*) \dasharrow G(3,k^5)$}

Let $r$ be an integer and let $F$ be a vector
space of dimension $2r+1$. An element $\omega$ in $\wedge^2 F^*$
yields a skew-symmetric map $F \rightarrow F^*$ which will be denoted
$L_\omega$. The rank, image, and kernel of $\omega$ will be those of
$L_\omega$. If $f_1,f_2 \in F$, $\omega(f_1,f_2)$ will denote the
number $L_\omega(f_1)(f_2)$.

\begin{lemm}
Let $\omega \in \wedge^2 F^*$ of rank $2t$ and $U \subset F$ a linear
subspace of dimension $2r+1-t$ and such that $\omega \bot \wedge^2 U$. Then
\begin{itemize}
\item
If $u \in U$, then $\omega(u,U) \equiv 0$.
\item
$\ker \omega \subset U$.
\end{itemize}
\end{lemm}
\pr
Taking a basis of $F$ containing a basis of $U$ and decomposing
$\omega$ along this basis, one checks that the condition 
$\omega \bot \wedge^2 U$ is equivalent to 
$\forall u,v \in U, \omega(u,v)=0$, proving the first point.

Therefore, we have $L_\omega(U) \subset U^\bot$, and since
$2t = \mathtt{rg}(L_\omega) \leq \mathtt{rg}({L_\omega}_{|U}) + t$,
we have $\mathtt{rg}({L_\omega}_{|U}) = t$ and so
$L_\omega(U)=U^\bot$. Since moreover
$L_\omega$ is skew-symmetric,
it follows that
$\ker L_\omega = {(\im\ L_\omega)}^\bot \subset {(U^\bot)}^\bot = U$.
\qed

\lpara

\begin{nota}
Let $\omega_1,\omega_2 \in \wedge^2 F^*$. We denote
$$
l(\omega_1,\omega_2) := \{ f \in F : \forall u \in \ker \omega_1,
\omega_2(u,f)=0 \}.
$$
\end{nota}

\begin{lemm}
Assume $2(2r+1) = 5t$. Let
$\omega_1,\omega_2 \in \wedge^2 F^*$ with rank $2t$ be such that
\begin{enumerate}
\item
$\ker \omega_1 \cap \ker \omega_2 = \{0\}$, and
\item
$L_{\omega_2}(\ker \omega_1) \cap L_{\omega_1}(\ker \omega_2) = \{ 0 \}$.
\end{enumerate}
If a linear subspace $U \subset F$ of dimension $2r+1-t$ is such that 
$\wedge^2 U \bot \omega_i,i=1,2$, then 
$U = l(\omega_1,\omega_2) \cap l(\omega_2,\omega_1)$.
\label{conditions}
\end{lemm}
\noindent
We will see in lemma \ref{existence} that for the minimal possible
values of $r,t$, which are those of interest to describe Mukai's flop, 
$U = l(\omega_1,\omega_2) \cap l(\omega_2,\omega_1)$ satisfies indeed 
$\wedge^2 U \bot \omega_i,i=1,2$; this is not the case in general.\\
\pr
Let $u \in \ker \omega_1$.
By the previous lemma, we have $\ker \omega_1 \subset U$, so $u \in U$.
If $f \in U$, it follows that $\omega_2(u,f)=0$, so 
$f \in l(\omega_1,\omega_2)$ and 
$U \subset l(\omega_1,\omega_2)$. By symmetry,
we have also $U \subset l(\omega_2,\omega_1)$. By condition (1),
$$\dim \ker \omega_1
= \dim L_{\omega_2}(\ker \omega_1) 
= \dim L_{\omega_1}(\ker \omega_2) = 2r+1-2t,$$ 
and by condition (2),
$$ \dim  (L_{\omega_2}(\ker \omega_1) +  
\dim L_{\omega_1}(\ker \omega_2)) = 2(2r+1-2t).$$
Since we know that $U$ is orthogonal to this space, and since by
hypothesis $2r+1 - 2 (2r+1-2t) = 4t-(2r+1) = 2r+1 - t = \dim U$,
$U$ is exactly
the orthogonal of this space, proving the lemma.
\qed

\begin{nota}
Denote $U(\omega_1,\omega_2) :=  l(\omega_1,\omega_2) 
\cap l(\omega_2,\omega_1)$.
\label{U}
\end{nota}

\begin{lemm}
Assume $r=t=2$ and let $\omega_1,\omega_2 \in \wedge^2 F^*$ be arbitrary.
Then there exists $U \subset F$ of dimension 3 such that 
$\wedge^2 U \bot \omega_1,\omega_2$. Therefore, if
the two conditions of lemma \ref{conditions} are satisfied, then
$\wedge^2 U(\omega_1,\omega_2) \bot \omega_1,\omega_2$.
If moreover
$\omega_1',\omega_2'$ are linear combinations of $\omega_1,\omega_2$
which also satisfy the two conditions, then
$U(\omega_1',\omega_2')=U(\omega_1,\omega_2)$.
\label{existence}
\end{lemm}
\pr
The second claim is a consequence of the first and the lemma
\ref{conditions}. The third claim follows from the second since 
$\wedge^2 U(\omega_1,\omega_2) \bot \omega_1',\omega_2'$.
It is therefore enough to prove the first claim. 

Let 
$G = \wedge^2 F^* \oplus \wedge^2 F^* \simeq Hom(k^2,\wedge^2 F^*)$. 
There is a natural 
$GL_2 \times GL(F)$ action on
$G$. Let $G(3,F)$ denote the grassmannian of 3-spaces in $F$ and
consider the
incidence variety $I \subset G(3,F) \times \p G$
defined by $(U,[\omega_1,\omega_2]) \in I$ \iff 
$\wedge^2 U \bot \omega_1,\omega_2$. It is a closed projective
$GL_2 \times GL(F)$-stable variety. Therefore, its projection 
$p_2(I) \subset \p G$ also.

Now, let $f_1,\ldots,f_5$ be a basis of $F$ and
$f_1^*,\ldots,f_5^*$ be the dual basis of $F^*$. Set
$\omega_1 = f_4^* \wedge f_1^* + f_5^* \wedge f_2^*$ and
$\omega_2 = f_4^* \wedge f_2^* + f_5^* \wedge f_3^*$.
It is clear that if $U = \mathtt{Vect}(f_1,f_2,f_3)$, then
$\wedge^2 U \bot \omega_1,\omega_2$; therefore 
$[\omega_1,\omega_2] \in p_2(I)$. It is proved in 
\cite[proof of proposition 13 p.94]{kimura} that the 
$GL_2 \times GL(F)$-orbit through $[\omega_1,\omega_2]$ is dense
(it also follows from lemma \ref{orbites});
therefore $p_2(I)=G$ and the existence claim of the lemma is proved.
\qed

%******************************************************************************

\subsectionplus{$GL_2 \times GL_5$-orbits}

Let as above $F$ a 5-dimensional vector space over $k$.
In this subsection, I describe the 
$GL_2 \times GL(F)$-orbits in $Hom(k^2,\wedge^2F^*)$, and prove
where the previous rational map 
$U : \p Hom(k^2,\wedge^2F^*) \dasharrow G(3,F)$, defined
on the open orbit by notation \ref{U}, extends.

We start with a result of co-diagonalisation of 2-forms of maximal
rank~:
\begin{lemm}
Let $\omega_1,\omega_2 \in \wedge^2 F^*$
be forms such that 
$\forall (\alpha_1,\alpha_2) \in k^2 - \{(0,0)\}$, 
$\alpha_1 \omega_1 + \alpha_2 \omega_2$ has rank 4. Then there exists
a basis $f_1^*,\ldots,f_5^*$ of $F^*$ such that
$$
\begin{array}{rcl}
\omega_1 & = & f_2^* \wedge f_4^* + f_3^* \wedge f_5^*\\
\omega_2 & = & f_1^* \wedge f_5^* + f_3^* \wedge f_4^*.
\end{array}
$$
\label{codiagonalisation}
\end{lemm}
\pr
For $i \in \{1,2\}$ and $u,v \in F$, we denote
$\scal{u,v}_i := L_{\omega_i}(u)(v)$.

Assume that $\ker \omega_1 = \ker \omega_2$. Denote $K$ this
1-dimensional vector space. Then
$\omega_1,\omega_2$ belong to $\wedge^2 (F/K)^*$. The variety of
degenerate 2-forms in $(F/K)^*$ is a hypersurface, so there exists 
$(\alpha_1,\alpha_2) \in k^2 - \{(0,0)\}$ such that
$\alpha_1 \omega_1 + \alpha_2 \omega_2$ is degenerate, contradicting
the hypothesis of the lemma.

We consider $0 \not = f_1\in \ker \omega_1$
and $0 \not = f_2 \in \ker \omega_2$; $f_1$ and $f_2$ are 
therefore not colinear.

\lpara

Assume now that $L_{\omega_1}(f_2)$ and $L_{\omega_2}(f_1)$ are
colinear. Denote $I$ this common image. The map 
$F^* \rightarrow F^*/I$ induces a map
$\wedge^2 F^* \rightarrow \wedge^2 (F^*/I)$
Let
$\overline \omega_i \in \wedge^2 {(I^\bot)}^*$ 
denote the image of $\omega_i$ under this map.
Both $\overline \omega_1$ and $\overline \omega_2$ vanish on
$f_1,f_2$. Therefore, they are proportional 2-forms~: let
$\alpha_1 \overline \omega_1 + \alpha_2 \overline \omega_2$ be a
non-trivial relation. Since $I^\bot$ is an isotropic subspace for
$\alpha_1 \omega_1 + \alpha_2 \omega_2$, this form does not have rank
4, contradicting the hypothesis.

We set $f_5^* = L_{\omega_2}(f_1)$ and
$f_4^* = L_{\omega_1}(f_2)$; $f_4^*$ and $f_5^*$ are 
therefore not colinear. Note that
$\scal{f_4^*,f_1} = \scal{f_2,f_1}_1 = 0$ because
$L_{\omega_1}(f_1)=0$, and that
$\scal{f_4^*,f_2} = \scal{f_2,f_2}_1 = 0$; therefore, 
$f_4^* \in \scal{f_1,f_2}^\bot$, and similarly
$f_5^* \in \scal{f_1,f_2}^\bot$.

\lpara

We now let $[\omega_i]$ be the composition
$F \stackrel{L_{\omega_i}}{\rightarrow} F^*
\rightarrow F^*/\scal{f_4^*,f_5^*}$. I claim that
$\ker [\omega_i] = \scal{f_4^*,f_5^*}^\bot$.

Note that $\im L_{\omega_i} = f_i^\bot \supset \scal{f_4^*,f_5^*}$.
I will prove the claim when $i=1$. Both spaces are
3-dimensional and contain $\scal{f_1,f_2}$. So let $f$ such that
$L_{\omega_1}(f) = f_5^*$, and let us see that 
$f \in \scal{f_4^*,f_5^*}^\bot$. Since $L_{\omega_1}(f) = f_5^*$ by
assumption, $\scal{f_5^*,f} = 0$. Similarly,
$\scal{f_4^*,f} = \scal{f_2,f}_1 = - \scal{f_5^*,f_2}
= - \scal{f_1,f_2}_2 = 0$, since $L_{\omega_2}(f_2) = 0$. So the claim
is proved.

Looking at the surjective maps
$\scal{f_4^*,f_5^*}^\bot 
\stackrel{L_{\omega_1},L_{\omega_2}}{\longrightarrow}
\scal{f_4^*,f_5^*}$, one proves that there exists
$f_3 \in \scal{f_4^*,f_5^*}^\bot$ such that
$L_{\omega_1}(f_3) \in \scal{f_5^*} - \{0\}$ and 
$L_{\omega_2}(f_3) \in \scal{f_4^*} - \{0\}$. 
Up to scaling $f_1$ (and so $f_5^* = L_{\omega_1}(f_1)$) and $f_2$
(and so $f_4^*$), we may assume that
$L_{\omega_1}(f_3) = f_5^*$ and 
$L_{\omega_2}(f_3) = f_4^*$.

\lpara

Up to now, the vectors $f_1,f_2,f_3,f_4^*,f_5^*$ were determined, up
to a scale, by $\omega_1$ and $\omega_2$. We now make a more significant
choice for $f_4$~: let $f_4 \in {(f_5^*)}^\bot$ such that
$\scal{f_4^*,f_4} = 1$. Note that this implies
$\scal{f_3,f_4}_2 = \scal{f_2,f_4}_1 = 1$, by definition of $f_3$ and
$f_4^*$. We moreover choose $f_5 \in {(f_4^*)}^\bot$ such that
$\scal{f_5^*,f_5} = 1$ and
$\scal{f_4,f_5}_1 = \scal{f_4,f_5}_2 = 0$. Note that this implies
$\scal{f_5,f_3}_1 = \scal{f_5,f_1}_2 = -1$.

It is then easy to check that for $i \in \{1,\ldots,5\}$, we have 
$f_5^*(f_i) = \delta_{i,5}$ and
$f_4^*(f_i) = \delta_{i,4}$. So it will not conflict notations to
consider the dual basis $(f_1^*,\ldots,f_5^*)$ of the basis
$(f_1,\ldots,f_5)$ of $F$. In this dual basis, $\omega_1$ and
$\omega_2$ are as in the proposition.
\qed

\para

Let $(f_1^*,\ldots,f_5^*)$ be a basis of $F^*$.
Let $\omega_1,\omega_2$ denote the forms
$f_2^* \wedge f_4^* + f_3^* \wedge f_5^*,
f_1^* \wedge f_5^* + f_3^* \wedge f_4^*$.
We now classify the $GL_2 \times GL(F)$-orbits in
$Hom(k^2,\wedge^2 F^*)$.

\begin{lemm}
There are eight $GL_2 \times GL_5$-orbits in $Hom(k^2,\wedge^2 F^*)$. The
following array gives elements in each orbit, its dimension and a
label.
$$
\begin{array}{|c|c|c|c|}
\hline
\mbox{label} & f((1,0)) & f((0,1)) & \dim \\
\hline
A_2 + 2A_1 & \omega_1 & \omega_2 & 20 \\
\hline
A_2 + A_1 & \omega_1 & f_1^* \wedge f_2^* & 18 \\
\hline
A_2      & \omega_1 & f_2^* \wedge f_4^* & 16 \\
\hline
3A_{1,a} & \omega_1 & f_2^* \wedge f_3^* & 15 \\
\hline
3A_{1,b} & f_1^* \wedge f_2^* & f_1^* \wedge f_3^* & 12\\
\hline
3A_{1,c} & \omega_1 & \omega_1 & 11\\
\hline
2A_1     & f_1^* \wedge f_2^* & f_1^* \wedge f_2^* & 8\\
\hline
A_1      & 0 & 0 & 0\\
\hline
\end{array}
$$
Finally, the closure of an orbit $\cal O$ contains the orbit ${\cal O}'$
\iff $\cal O$ lies above ${\cal O}'$ in this array, except that the
closure of the orbit labelled $3A_{1,b}$ does not contain the orbit
labelled $3A_{1,c}$. 
\label{orbites}
\end{lemm}
\pr
Granting the classification of the orbits, I leave it to the reader to
check the dimensions of the orbits and the decomposition of their closures.

So let again $F = k^5$ and $f \in Hom(k^2,\wedge^2 F^*)$. If the rank
of $f$, as a morphism of vector spaces, is one, then there are three
cases (labelled $A_1,2A_1,3A_{1,c}$), according to the rank 
(as an element of $\wedge^2 F^*$) of a generic element of its
image.

Assume $f$ has rank two.
If all non-vanishing elements of the image of $f$ have rank 4, then,
by lemma \ref{codiagonalisation}, we are in case $A_2+2A_1$. If all
these elements are degenerate, then it is well-known that we are in
case $3A_{1,b}$.

Otherwise, we may assume that $f((1,0))$ has rank 4 and
$\omega := f((0,1))$ has
rank 2. There is a basis $f_1,\ldots,f_5$ of $F$ such that in terms
of the dual basis $f_1^*,\ldots,f_5^*$, $f((1,0)) = \omega_1$. The
kernel of $L_{\omega_1}$ is generated by $f_1$. Consider the
4-dimensional subspace $F':=\scal{f_2^*,f_3^*,f_4^*,f_5^*} \subset F^*$
and the 5-dimensional projective space 
$\p \wedge^2 F'$ containing the 4-dimensional quadric $G(2,F')$ of
classes of elements of rank 2. The generic element 
$\omega_1 \in \wedge^2 F'$ defines a polar hyperplane 
(with respect to the quadric)
in $\p \wedge^2 F'$, which will be denoted
$H$. Note that $H \cap G(2,F')$ is a smooth 3-dimensional quadric.
\begin{itemize}
\item
Assume first that $L_\omega$ does not vanish on $f_1$ and let $g_1^*$
be an element in $\im\ L_\omega$ not vanishing on $f_1$.
Let $g_2^* \not = 0$ be an element in 
$\im\ L_\omega \cap \im\ L_{\omega_1}$.
The variety of classes of elements of the form 
$[g_2^* \wedge g^*],g^* \in F'$ is a $\p^2$ in the quadric
$G(2,F')$. Therefore, it can't be included in $H$. Let $g_4^* \in F'$ such
that $[g_2^* \wedge g_4^*] \not \in H$; the projective line through
$[g_2^* \wedge g_4^*]$ and $[\omega_1]$ is therefore a secant line~:
let $[g_3^* \wedge g_5^*]$ be an element in the intersection of this
line and $G(2,F')$. We can assume
$\omega_1 = g_2^* \wedge g_4^* + g_3^* \wedge g_5^*$ and
$\omega = g_1^* \wedge g_2^*$; therefore, we are in the case labelled
$A_2+A_1$.
\item
Assume now that $L_\omega(f_1) = 0$. In this case, both $\omega$ and
$\omega_1$ are in $\wedge^2 F'$. There are two $GL(F')$-orbits for the
projective line through $\omega$ and $\omega_1$~: either it is a
secant line to the quadric $G(2,F')$, either it is a tangent line;
this corresponds to the cases $A_2$ and $3A_{1,a}$.
\vspace{-.65cm}
\end{itemize}
\qed

\para

Recall the rational map $U$ of notation \ref{U}. It is a model for the
restriction of the Mukai flop of the second kind to a tangent space,
so it is interesting
to know where it is defined.

\begin{lemm}
The open orbit is the locus where $U$ is defined.
\end{lemm}
\pr
Let as before $\omega_1 = f_2^* \wedge f_4^* + f_3^* \wedge f_5^*$,
and let
$\omega_2(t) = f_1^* \wedge f_5^* + t.f_2^* \wedge f_3^*$. Let 
$f : k^2 \rightarrow \wedge^2 F^*$ be defined by
$f((1,0)) = \omega_1$ and $f((0,1)) = \omega_2(t)$; we have
$U(f) = \scal{f_5^*,f_2^*}$. 

The same construction with 
$\omega'_2(t) = f_1^* \wedge f_5^* + t.f_3^* \wedge f_4^*$ yields
$U(f') = \scal{f_5^*,f_4^*}$. Now, since $\omega_2(t)$ and
$\omega_2'(t)$ converge to $f_1^* \wedge f_5^*$, this proves that $U$
is not defined at the point $f_0$ defined by
$f_0((1,0)) = f_2^* \wedge f_4^* + f_3^* \wedge f_5^*$ and
$f_0((0,1)) = f_1^* \wedge f_5^*$.

Therefore, the indeterminacy locus of $U$ contains the orbit labelled
$A_2+A_1$; since it is closed, it contains all the orbits but the open
one.
\qed

%******************************************************************************

\subsectionplus{Unicity of the equivariant rational map}

Recall that $F$ is a 5-dimensional vector space over $k$.
In this subsection, I show that the rational map 
$U : Hom(k^2,\wedge^2 F^*) \dasharrow G(3,F)$ of notation
\ref{U} is the unique $(GL_2 \times GL(F))$-equivariant rational map
$Hom(k^2,\wedge^2 F^*) \dasharrow G(3,F)$. This is a result analogous
to proposition \ref{unicite_spin}, and the strategy of proof is the
same~: we caracterize its graph as an orbit of minimal dimension.

Let $O$ denote the open $(GL_2 \times GL(F))$-orbit in
$Hom(k^2,\wedge^2 F^*)$.

\begin{lemm}
In $O \times G(3,F)$, there is a unique 20-dimensional orbit.
\label{unique_dim20}
\end{lemm}
\pr
The set of $(f,\alpha)$ with $\alpha = U(f)$ is such an orbit. Let
$(f,\alpha)$ be in an orbit of dimension 20. Since $O$ is homogeneous,
we can assume that 
$f((1,0)) = \omega_1 = f_2^* \wedge f_4^* + f_3^* \wedge f_5^*$ and
$f((0,1)) = \omega_2 = f_1^* \wedge f_5^* + f_3^* \wedge f_4^*$. Let
$G_0$ denote the stabilizor of $f$ in $GL_2 \times GL(F)$; we have
$G_0 = \{1\} \times G_1$, with
$G_1 = \left \{ \left (
\begin{array}{ccccc}
t & 0 & 0 & 0 & 0 \\
0 & t & 0 & 0 & 0 \\
0 & 0 & t & 0 & 0 \\
a & c & d & t^{-1} & 0 \\
b & d & a & 0 & t^{-1}
\end{array}
\right )
: t \in k^*; a,b,c,d \in k \right \}
$
(these matrices express the action of $G_1$ on $F^*$ in the dual basis
$f_1^*,\ldots,f_5^*$). This fact can be proved by a direct
computation; we will only use the obvious fact that $G_1$ stabilizes
$f$.

Let $F_\alpha \subset F$ denote the 3-dimensional subspace
corresopnding to $\alpha$, and let $F_\alpha^\bot \subset F^*$ denote
its orthogonal.
Since $(f,\alpha)$ belongs to an orbit of dimension 20, $G_1$ must
stabilize $\alpha$. Assume there exists $f^* \in F_\alpha^\bot$, with
$f^* = \sum x_i f_i^*$ and $x_1 \not = 0$. The action of $G_1$ implies
that $f_5^*$ belongs to $F_\alpha^\bot$, and so also
$f_4^*$. Therefore we have a contradiction. The same contradiction
arises if $F_\alpha^\bot$ contains a form with non-vanishing
coefficient along $f_2^*$. From this it follows easily that
$F_\alpha^\bot$ is generated by $f_4^*$ and $f_5^*$, so 
$\alpha = U(f)$.
\qed

\begin{prop}
There is a unique $(GL_2 \times GL(F))$-equivariant rational map
$Hom(k^2,\wedge^2 F^*) \dasharrow G(3,F)$.
\label{unicite_sl5}
\end{prop}
\pr
$U$ is such a rational map. Let $u$ denote any
$(GL_2 \times GL(F))$-equivariant rational map
$Hom(k^2,\wedge^2 F^*) \dasharrow G(3,F)$.

Recall that $O \subset Hom(k^2,\wedge^2 F^*)$ denotes the open orbit;
since $u$ is equivariant, it is defined on $O$ and surjective. The
variety $\{(f,u(f)):f\in O\}$ is a 20-dimensional orbit in 
$O \times G(3,F)$; therefore, by lemma \ref{unique_dim20},
it is equal to the variety $\{(f,U(f)):f\in O\}$.
\qed

%******************************************************************************

%******************************************************************************

%******************************************************************************

\sectionplus{Tangency in Scorza varieties}

\label{scorza}

In a projective space, given a point $x$ and a non-vanishing
tangent vector 
$t \in T_xX$, there is a unique line $l$ through $x$ and such that 
$t \in T_xX$. Similarly, given a non-vanishing
cotangent form $f \in T^*_xX$, there
is a unique hyperplane $h$ such that $f$ vanishes on $T_xh$.
Therefore, a tangent vector defines a line and a cotangent form a 
hyperplane. This will be extended to
a projective space over a composition algebra in this section.
Both of these maps will be also defined using Jordan algebras.

\subsectionplus{Notations for Scorza varieties}

\label{def_scorza}

Let $\a$ be a composition algebra over $\C$, of dimension $a$. If $n$
is an integer, let $H_n(\a)$ denote the space of $(n \times n)$
hermitian matrices with coeficients in $\a$.

Let
$$\fonction{\nu_2}{\a^n}{H_n(\a)}{(z_1,\ldots,z_n)}
{(z_i\overline z_j)_{1 \leq i,j \leq n}}$$
be the map generalizing that of section \ref{spin10}.
Recall from \cite{projective} that in $H_n(\a)$ there is a notion of rank.
The variety of rank $n-1$-elements is a hypersurface; let $\det$ denote a
reduced equation of this hypersurface. The variety of rank one
matrices may be described, by
\cite[theorem 3.1 (4) and proposition 4.2]{projective}, as the closure
$X = \overline {\{ [\nu_2(1,z_2,\ldots,z_n)]\ :\ z_i \in \a \}}$.

\lpara

Scorza varieties were defined and classified by F. Zak as varieties
having some extremal properties with respect to their secant varieties
\cite{zak,scorza}. For our purpose, the following theorem
will serve as a definition~:
\begin{theo}[Zak]
Let $a \in \{1,2,4,8\}$ and $n$ be integers.
A Scorza variety of type $(n,a)$
is a pair $(V,X)$, where $V$ is a $\C$-vector space,
and $X \subset \p V$ is a projective variety projectively
isomorphic to the variety of classes of rank one matrices in the
projectivisation of the space $H_n(\a)$ (with $\dim \a = a$).
\end{theo}

$X$ is a kind of projective space; moreover, one can define a dual
``projective space'' $X^\dual \subset \p V^*$,
non-canonically isomorphic with $X$,
and an incidence relation for
$(x,h) \in X \times X^\dual$ denoted $x \vdash h$. 
In fact, $X^\dual$ is the variety of hyperplanes containing $n-1$
general tangent spaces to $X$ and $x \vdash h$ \iff $T_xX \subset h$.
For $h \in X^\dual$,
the Schubert cell of $x$'s incident to $h$ will be denoted $C_h$.
The quadratic representation corresponding to the Scorza variety
$(V^*,X^\dual)$ will be denoted $U^\dual$; therefore, $U^\dual$ is a
quadratic map $V^* \rightarrow Hom(V,V^*)$.

For the convenience of the reader, I recall, given $a$ and $n$, the
corresponding Scorza varieties and their automorphism group ($G(n_1,n_2)$ 
denotes the grassmannian of $n_1$-dimensional subspaces in $\C^{n_2}$.

$$
\begin{array}{|c|c|c|c|c|c|}
\hline
a & \a & V & X & X^\dual & Aut(X)\\
\hline
1 & \C & S^2\C^n & \p^{n-1} & {(\p^{n-1})}^\dual & PGL_n\\
\hline
2 & \C \oplus \C & \C^n \otimes \C^n & \p^{n-1} \times \p^{n-1} &
{(\p^{n-1})}^\dual \times {(\p^{n-1})}^\dual &
PGL_n \times PGL_n\\
\hline
4 & M_2(\C) & \wedge^2 \C^{2n} & G(2,2n) & G(2n-2,2n) & PGL_{2n}\\
\hline
8 & \oc & \dim 27 & E_6/P_1 & E_6/P_6 & E_6\\
\hline
\end{array}
$$

In each case, it is well-known that there is a Mukai flop
$T^*X \dasharrow T^*X^\dual$. One aim of the rest of the article is
to describe this flop.

\lpara

Let $(X,V)$ be a Scorza variety of type $(n,a)$.
Recall from \cite[section 1]{hermitian} 
the ``quadratic representation''~: it is
a quadratic map $V \rightarrow Hom(V^*,V)$, canonically defined using only
$X$. If $A \in V$, we will denote $U_A \in Hom(V^*,V)$ the image of
$A$ under the quadratic representation.

In concrete terms, when we will have to compute a quadratic
representation in $V$, we will allways do the following. First, we
will identify $V$ with $H_n(\a)$. Second, we will choose the scalar
product $(A,B) = \trace(AB)$, which identifies $V$ and $V^*$. These
two choices will not affect the final result. Then, to
compute $U_A(B)$, for $A \in V$ and $B \in V^* \simeq V$, we will
allways manage to be in the situation when all the coefficients of $A$
and $B$ belong to an associative subalgebra of $\a$ (this holds, for
example, if $\a$ itself is associative). Then we use the fact that
$U_A(B)$ is $ABA$, where juxtaposition stands for the usual product of
matrices \cite{hermitian}.

Recall also that for any integer $r<n$ there is a well-defined variety 
$G_\a(r,X)$
parametrizing Scorza subvarieties of type $(r,a)$ in $X$. To an element
$A \in \p V$ of rank $r$ is associated a subvariety
$X_A \in G_\a(r,X)$ and its linear span in $\p V$ is denoted $\Sigma_A$
\cite[proposition 1.3]{hermitian}.

\lpara

As explained in \cite{projective} and \cite{hermitian}, the Scorza
varieties admit a model over $\Z$, and the quadratic representation is
defined over $\Z$. Therefore, all the following constructions are valid
on this base, and we get a description of Mukai flops over $\Z$. For
the clarity of redaction, I will work over $\C$, since it is the usual
context of Mukai flops.
\lpara

In the following, $(V,X)$ will be a Scorza variety of type $(n,a)$,
and $G$ denotes the automorphism group of $X$.

\subsectionplus{A generic tangent vector defines a line}

Let $x \in X$ and let $L_x \subset V$ be the line it represents. 
We have $T_xX = Hom(L_x,\widehat {T_xX}/L_x)$.
Let $t \in T_xX$;
in the next
proposition, I say that $T \in \widehat{T_xX}$ represents $t$ if the
morphism $t \in Hom(L_x,\widehat {T_xX}/L_x)$ has image the line
generated by the class modulo $L_x$ of $T$.

By \cite[proposition 1.5]{hermitian}, the $\a$-lines through a point 
$x \in X$ are naturally parametrized by a subvariety of 
$\p (V/\widehat {T_xX})$. I say that a representative of an $\a$-line
through $x$ is $L$ (with $L \in V/\widehat {T_xX}$) if the class of $L$
in $\p (V/\widehat {T_xX})$ corresponds to $l$.

\begin{theo}
Let $x \in X$ and $t \in T_xX$ generic. There exists a unique
$\a$-line $l \in G_\a(2,X)$ such that $x \in l$ and $t \in T_xl$.
A representative for $l$ in $V/\widehat {T_xX}$ is
$L = [U_T(A)]$, if 
$T \in \widehat{T_xX}$ represents $t$ and $A$ is a generic element in
$V^*$.
\label{tangent}
\end{theo}
\begin{nota}
Let $\nu_x^+$ denote the quadratic map 
$T_xX \rightarrow V/\widehat{T_xX},T \mapsto U_T(A)$ of this theorem.
\end{nota}
\pr
Let $x \in X$ and $t \in T_xX$ be generic. Let $T$ represent $t$.
Then $T$ has rank two, so by \cite[proposition 1.4]{hermitian},
$T$ defines the $\a$-line $X_T$. We will prove that $X_T$ is the
unique $\a$-line with the properties of the proposition. To this end,
we assume that $n=3$ to simplify notations, since larger values of $n$
would not change the argument.

We assume that $V = H_n(\a)$ and $X$ is the variety of rank one matrices.
By \cite[proposition 1.3]{hermitian}, we can
assume $x=\matttr 100000000$ and $T=\matttr 110100000$. Then 
$X_T$ is the set of rank one
matrices of the form 
$\mattt**0**0000$. We therefore check that $x \in X_T$ and 
$T \in \widehat{T_xX_T} = \matttr **0 *00 000$.

Conversely, let $l \in G_\a(2,X)$ such that $x \in l$ and $t \in
T_xl$. Let $B \in \p V$ of rank 2 such that $l = X_B$.
Then since $T$ represents $t$ and $t \in T_xl$, we have 
$T \in \widehat{T_xl} \subset \Sigma_B$.
By \cite[proposition 1.4]{hermitian},
$\Sigma_B = \Sigma_T$ and $l=X_T$.

The fact that $L = [U_T(A)] \in V/\widehat {T_xX}$
is a representative for $l$ follows
from the fact that by \cite[proposition 1.4]{hermitian} again,
$\Sigma_T$ is the image of $U_T$, and the fact that
the isomorphism of 
\cite[proposition 1.5]{hermitian} maps the $\a$-line $X_T$ on the
line $\im\ U_T/\widehat{T_xX} \subset V/\widehat{T_xX}$.
\qed

\lpara

Let $d$ be an integer and $\a$ a composition algebra; recall the map
$\nu_2:\a^d \rightarrow H_d(\a)$ defined in subsection
\ref{def_scorza}. 
Its projectivisation $\overline \nu_2 : \a^d \dasharrow \p H_d(\a)$
may be considered as a kind of quotient map
$\a^d \dasharrow \p^{d-1}_\a$ \cite[subsection 3.4]{projective}.

\begin{coro}
There are identifications of $T_xX$ with $\a^{n-1}$
and $V/\widehat {T_xX}$ with $H_{n-1}(\a)$ such that $\nu_x^+$
identifies with $\nu_2:\a^{n-1} \rightarrow H_{n-1}(\a)$.
\end{coro}
\pr
With the notations of the previous proof,
to see that $\nu_x^+$ identifies with $\nu_2$, we choose the scalar
product $(A,B) \mapsto \trace(AB)$ on $V=H_n(\a)$, which identifies $V$ and
$V^*$, and moreover we choose $A \in V^*$ to be the linear form
corresponding to the identity matrix in $V=H_n(\a)$. 
Then, by subsection \ref{def_scorza}, if
$T=\matttr t{\overline z_1}{\overline z_2} {z_1}00 {z_2}00$, then
$U_T(A)=T^2=\matttr *** *{N(z_1)}{z_1\overline z_2} 
*{z_2\overline z_1}{N(z_2)}$. 
Therefore, $\nu_x^+$ identifies with $\nu_2$.
\qed

\subsectionplus{A generic cotangent form defines a hyperplane}

A cotangent form $f \in T_xX^*$ is an element
$f \in Hom(L_x^*,(\widehat{T_xX}/L_x)^*)$. 
I say that $\tilde f \in V^*$ represents
$f$ if $\tilde f_{|\widehat{T_xX}}$ generates the image of $f$.
Recall (subsection \ref{def_scorza}) that for $h \in X^\dual$, $C_h$
denotes the Schubert cell in $X$ defined by $h$. Let
$\mu:T^*X \dasharrow T^*X^\dual$ denote the Mukai flop and 
$\pi:T^*X^\dual \rightarrow X^\dual$ the projection.

\begin{theo}
Let $x \in X$, $x_0 \in L_x - \{0\}$, and $f \in T_x^*X$ generic. 
There exists a unique 
$h \in X^\dual$
such that $f$ vanishes on $T_xC_h$. If $\tilde f \in V^*$ 
represents $f$, then a representative of $h$ is 
$U^\dual_{\tilde f}(x_0) \in (V/\widehat{T_xX})^*$.
Finally, $\pi \circ \mu(x,f)=h$.
\label{cotangent}
\end{theo}
\begin{nota}
Let $\nu_x^-$ denote the quadratic map 
$T_x^*X \rightarrow (V/\widehat{T_xX})^*,
\tilde f \mapsto U^\dual_{\tilde f}(x_0)$ of this theorem.
\end{nota}
\pr
The last claim follows from the first and \cite{dual},
where it is proved that 
$\pi \circ \mu(x,f)$ is the only $h \in X^\dual$ such that $f$
vanishes on $T_xC_h$.

To simplify notations, we assume in the proof that
$V = H_3(\a)$ and we identify $V$ and $V^*$ via the scalar product
$(A,B)=\trace(AB)$.
Assume as before that $x = \mattt 100 000 000$ and that
$f \left ( \mattt t{\overline z_1}{\overline z_2} {z_1}00 {z_2}00
\right ) = \re(z_2)$.

Let $h_0 = \mattt 000 000 001$. Let $y = \mattt tz0 {\overline z}u0 000$,
with $t,u \in \C$ and $z \in \a$ and $tu-N(a)=0$, be an element of
$\widehat X$. By \cite[theorem 3.1 (4) and proposition 4.2]{projective},
$\widehat X = \overline {\{ [\nu_2(1,z_1,z_2)] : \ z_i \in \a \}}$.
We deduce that if $(m_{i,j}) \in \widehat X$, then the minors
$m_{i,i}m_{j,j}-N(m_{i,j})$ vanish. It follows that if $t \not = 0$,
$T_y \widehat Y$ is orthogonal to $h_0$.

Therefore, by continuity, the
intersection of $\mattt **0 **0 000$ with $X$ lies in the Schubert
cell $C_{h_0}$, and for dimension reasons we have equality. This shows
that $f$ vanishes on $T_x C_{h_0}$. Therefore $h=h_0$.

Finally, let
$\tilde f=\mattt 0{\overline z_1}{\overline z_2} {z_1}t{\overline z} 
{z_2}zu$ be a linear form
($t,u\in \C$ and $z,z_1,z_2 \in \a$ are arbitrary); then 
$U^\dual(\tilde f).x = \tilde f.x.\tilde f 
= \tilde f \matttr 0{\overline z_1}{\overline z_2} 000 000 =
\matttr 000 0{N(z_1)}{z_1\overline z_2} 0{z_2\overline z_1}{N(z_2)}$,
so that $\nu_x^-$ identifies with $\nu_2$. Moreover, if $z_1=0$ and
$z_2=1$, then $\tilde f$ represents $f$ and we have 
$[\nu_x(\tilde f)] = h_0$, as claimed.
\qed

\lpara

In the proof of the theorem, we showed~:

\begin{coro}
There are identifications of $T_xX^*$ and
$(V/\widehat{T_xX})^*$ with $\a^{n-1}$ and $H_{n-1}(\a)$ such that
$\nu_x^-$ identifies with $\nu_2:\a^{n-1} \rightarrow H_{n-1}(\a)$.
\end{coro}

\subsectionplus{The variety of lines through a point in
$\p^n_\a$ as Fano variety of maximal linear subspaces of 
$\p^{n-1}_\a$}

\label{variety}

The goals of this subsection are propositions \ref{fano_tangent} and
\ref{fano_cotangent}.
\lpara

The normal bundle to $X$ in $\p V$ twisted by $(-1)$ will be denoted
$N$ and let $\pi: N \rightarrow X$ 
(resp. $\overline \pi : \p N \rightarrow X$) be the structure map 
of this vector
bundle (resp. its projectivisation). Similarly, let $\psi$ and
$\overline \psi$ denote the natural maps $TX(-1) \rightarrow X$ and 
$\p TX(-1) \rightarrow X$. Let $x \in X$; the quotient map 
$V \rightarrow V/\widehat{T_xX} = N_x$ will be denoted $\pi_x$.
The normal bundle $N$ admits an interesting
subvariety~: the image of $\widehat X$. 
This variety will be denoted $N(X)$~: by
definition, the fiber $N(X)_x := \pi^{-1}(x) \cap N(X)$ is
$\pi_x(\widehat X)$. Recall \cite{hermitian} that $(N_x,\p N(X)_x)$ is
a Scorza variety of type $(n-1,a)$.

Assume $a>1$.
Let $F(0,1,X)$ denote the variety of couples $(x,l)$ where
$x \in \p V$, $l \subset \p V$ is a projective line, and 
$x \in l \subset X$. The map which
sends a pair $(x,l) \in F(0,1,X)$ to 
$(x,t) \in \p TX$, where $t \in \p T_xX$ is the
projectivisation of the tangent vector of $l$ at $x$ shows that
$F(0,1,X)$ can be considered as a subvariety of $\p TX$. By 
\cite[lemma 1.2 and proposition 1.3]{hermitian}, $F(0,1,X)$ is
homogeneous.

The first interesting point is that 
$(\overline \psi^*N)_{|F(0,1,X)}$ admits a
subbundle included in $\overline \psi^{-1}(N(X))$.
For $(x,l) \in F(0,1,X)$ and
$x \not = y \in l$, define $T_y := \pi_x(\widehat {T_yX})$.

\begin{prop}
$T_y$ does not depend on $y \in l$ and
$T_y \rightarrow F(0,1,X)$ defines a rank $(ra/2 + 1)$-subbundle
of $(\overline \psi^*N)_{F(0,1,X)}$, entirely included in
$\overline \psi^{-1}(N(X))$.
\label{t_y}
\end{prop}
\pr
Assume for the simplicity of notations that $n=3$.
I use the fact 
that if $z_1,z_2,z_3$ generate an associative subalgebra of $\a$, then
$\nu_2(z_1,z_2,z_3) \in \widehat X$ \cite[proposition 4.2]{projective}.
The condition on $z_1,z_2,z_3$ holds for example if $\a$ 
itself is associative or if $z_1=1$, since in $\oc$, the
subalgebra generated by two elements is allways associative.

Let $x=\matttr100000000$ and $y=\matttr1{\overline z}0z00000$ in 
$H_3(\a)$, with $z \in \a$ and $N(z)=0$. Then the line through $[x]$
and $[y]$ in $\p H_3(\a)$ lies in $X$, because 
$x+ty=\nu_2(1,tz,0)$. Moreover, differentiating $\nu_2$, we have
$\widehat{T_xX} = \left \{ \mattt ****00*00 \right \}$ and
$\widehat{T_yX} = \left \{ \mattt {\re(u)}{u \overline z + \overline v}
{\overline w}
{z \overline u +v}{\re(zv)}{z\overline w}
{w}{w\overline z}{0} : u,v,w \in \a \right \}$.
It follows that 
$$\widehat {T_yX}/\widehat {T_xX} \simeq 
\left \{ \matdd{*}{z\overline w}{w\overline z}{0} : w \in \a \right \}.$$
Therefore, this space does not change if $y$ is replaced by a point of
the line through $x$ and $y$. Since $F(0,1,X)$ is homogeneous, this
holds for any of its elements. Therefore, $T_{(x,l)}$ is allways a 
$(a/2 +1)$-linear subspace of $V/\widehat {T_xX}$. It follows
that it is a subbundle of $(\overline \psi^*N)_{|F(0,1,X)}$, 
as it is locally the image of
the bundle $\widehat{T_yX}$ ($y$ a local section of $l$ different from
$x$) under a constant rank vector bundle map.

Moreover, since $\matttr0000*{z\overline w}0{w\overline z}0$
belongs to $\widehat X$, 
$\matdd *{z \overline w}{w\overline z}0$ belongs to $N(X)_x$, and so
the subbundle $T$ is included in $\overline \psi^{-1}(N(X)_x)$.
\qed

\lpara

Let $x \in X$ and denote $F(x,1,X)$ the variety of lines through $x$
and included in $X$.
Theorem \ref{tangent}
yields a quadratic map 
$\nu^+_x : \widehat {T_xX}/L_x \rightarrow N_x$, well-defined up
to a scale. Let
$\mu^+_x(.,.)$ denote its polarization.

\begin{prop}
The map 
$l \mapsto T_{(x,l)}$ defines a map between
$F(x,1,X)$ and
some components of the Fano variety of maximal linear
subspaces in $\p N(X)_x$. This map is an isomorphism when $a \not = 4$,
and is surjective with fibers isomorphic to $\p^1$ when $a=4$.
Moreover, let $0 \not = t \in T_xl$;
$T_{(x,l)}$ is the image of $\mu_x^+(t,.)$.
\label{fano_tangent}
\end{prop}
\noindent
In particular, in the case $\a = \oc$ this proposition proves that the
variety of lines in $X$ through a fixed point $x \in X$ is isomorphic
with a 10-dimensional spinor variety; this fact is proved in
\cite[prop 3.4 p.77]{manivel}, but I give here a direct proof which
makes it clear which isotropic spaces this spinor variety parametrizes.\\
\pr
First of all, by \cite[proposition 1.4]{hermitian}, there are two
$P$-orbits in $\p T_xX$, if $P$ denotes the stabilizor of $x$ in $G$.
Therefore, $F(x,1,X) = \p \{\nu_x = 0 \} \subset \p T_xX$. We know that
$\nu_x$ identifies with $\nu_2$ and in 
\cite[section 3.4]{projective}, I described the locus where $\nu_2$
vanishes; therefore, we get the following array ($G^+_Q(5,10)$ denotes
a 10-dimensional spinor variety and $Q$ an 8-dimensional 
projective quadric; $G(2,2n-2)$ is the grassmannian of 2-dimensional
spaces in $\C^{2n-2}$)~:

$$
\begin{array}{|c|c|c|c|}
\hline
V        & F(x,1,X) & \p N(X)_x & \frac a2(n-2)+1 \\
\hline
H_n(\cc) & \p^{n-2} \amalg \p^{n-2} & \p^{n-2} \times \p^{n-2} & n-1\\
\hline
H_n(\hc) & \p^1 \times \p^{2n-3} & G(2,2n-2) & 2n-3\\
\hline
H_3(\oc) & S^+ & Q & 5\\
\hline
\end{array}
$$
Let now $z \in \a$ such that $N(z)=0$; the image of 
$\mu^+_x((z,0,\ldots,0),.)$ is $$\left \{ \matddr *
{z\overline u_1,\ldots,z\overline u_{n-2}}
{
\begin{array}{c}
u_1 \overline z\\
\vdots\\
u_{n-2} \overline z
\end{array} }
0 : u_i \in \a\right \}\ \ ;$$ 
it is of dimension $1+\frac a2(n-2)$,
so it is a maximal linear subspace of $N(X)_x$. 
If $l \in F(x,1,X)$, the fact that the image of $\mu_x^+(l,.)$
is $T_{(x,l)}$ is a consequence of the formula for $\nu_x^+$ and the
computation of $T_{(x,l)}$ made in the proof of proposition \ref{t_y}.

In the case when $\a = \oc$, proposition \ref{s+-} shows that
the map of the proposition is an isomorphism. I leave it to the reader
to check that in case $\a = \cc$, it is an isomorphism, and in case 
$\a = \hc$, it has fibers isomorphic with $\p^1$.
\qed

\lpara

Let $\nu^-_x:(\widehat{T_xX}/L_x)^* \rightarrow N_xX^* \subset V^*$ 
be the quadratic map of theorem \ref{cotangent} and $\mu_x^-$ its
polarization. We know that
$(N_x,\p N(X)_x)$ is a Scora variety of type $(n-1,a)$; let
$\p N(X)_x^\dual \subset \p N_x^*$ denote its dual Scorza variety.

We have a similar result for the cotangent space~:

\begin{prop}
The map 
$l \mapsto \im\ \mu_x^-(t,.)$, where $0 \not = t \in T_xl$ 
defines a morphism between
$F(x,1,X)$ and
some components of the Fano variety of maximal linear subspaces in
$\p N(X)_x^\dual$. It is an isomorphism if $a \not = 4$, and is
surjective with
fibers isomorphic to $\p^1$ if $a=4$.
\label{fano_cotangent}
\end{prop}

\lpara

From the array in the proof of proposition \ref{fano_tangent}, we see
that in case $n=3$, $\p N(X)_x$ is a smooth quadric of dimension
$a$. So there are two families of maximal linear subspaces in
$N(X)_x$. In case $\a = \cc$, the two families are described by
proposition \ref{fano_tangent}. But in case $a \geq 4$, we only
get one family. The other family comes with
proposition \ref{fano_cotangent}, because we can
use the canonical isomorphism
$\p N(X)_x^\dual = \p N(X)_x$ which holds since $\p N(X)_x$ is a smooth
quadric. One can check that we indeed find two different families with
the two dual constructions of propositions \ref{fano_tangent} and
\ref{fano_cotangent}.

\subsectionplus{The tangent bundle to the variety of lines in a Severi variety}

\label{para_tangent_severi}

In this subsection, we prepare the description of the Mukai flop of the
second kind. Let $X \subset \p V$ be a Scorza scheme
of type $(3,a)$ (these schemes are called Severi varieties in
\cite{zak}) and assume $a \geq 2$. Note \cite{projective}
that if $a=2$, then $X$ is
isomorphic with $\p^2 \times \p^2 \subset \p^8$ and that if $a=4$,
then $X$ is isomorphic with the grassmannian $G(2,6) \subset \p^{14}$ 
of 2-dimensional subspaces of $\C^6$. 

Let $Y$ denote an irreducible component of the variety of projective lines in
$\p V$ which are included in $X$. If $a=2$, then 
$Y \simeq {(\p^2)}^\dual \times (\p^2)$ and if $a=4$, then $Y$
is isomorphic with the flag variety $F(1,3,6)$ of 1-dimensional subspaces
included in a 3-dimensional subspace included in a fixed $\C^6$. 
If $a=8$, then it follows from \cite[theorem 4.3 p.82]{manivel} 
that $Y$ is the
quotient $G/P_3$, where $G$ is a simply-connected group of type $E_6$
and $P_3$ is the parabolic subgroup corresponding to the simple root
$\alpha_3$. Therefore, the Mukai flop of the second kind is a rational
map $T^* Y \dasharrow T^* Y^\dual$, where $Y^\dual = G/P_5$.

The aim of this subsection is to describe the tangent bundle $TY$. As
before, this will be done in a unified way for all Severi varieties
with $a \geq 2$ (if $a=1$, the variety $Y$ is empty).

\para

Let us start with an easy lemma. Let $\det(.,.,.)$ be the polarization
of the degree 3 polynomial $\det$ (that is, the unique trilinear
symmetric form such that $\forall v \in V,\det(v,v,v)=6\det(v)$).
\begin{lemm}
Let $X$ be a Severi variety and $x \in X$. Then we have
$$\widehat{T_xX} = \{ v : \forall w \in V,\det(x,v,w)=0 \}.$$
\label{tangent_severi}
\vspace{-.6cm}
\end{lemm}
\pr
By \cite[propositions 3.5 and 4.2]{projective}, the ideal of $X$ is
generated by the quadratic equations $\det(x,x,.)=0$. Therefore, by
differentiation, we get the given equations for the tangent space at
$x$.
\qed

Now, let $\alpha \in Y$. The 2-dimensional linear space it represents
will be denoted $L_\alpha$. We set
$$
\begin{array}{rcl}
S_\alpha & := & \langle T_x \widehat X \rangle_{x \in L_\alpha - \{0\}}
\vspace{.1cm} \\
I_\alpha & := & \bigcap_{x \in L_\alpha - \{0\}} T_x\widehat X.
\end{array}
$$
It is clear that $S$ and $I$ are $G$-homogeneous 
subbundles of the
trivial bundle $V \otimes \co_Y$ over $Y$. We moreover consider the
quotient bundles defined by
$A_\alpha := I_\alpha / L_\alpha  ,  B_\alpha := S_\alpha / I_\alpha  ,
C_\alpha := V/S_\alpha$.

\begin{prop}
The ranks of the bundles $A,B,C$ are, respectively,
$3a/2-2,a+2,a/2+1$. There is a $G$-equivariant short exact sequence of
bundles
$$
0 \rightarrow Hom(L,A) \rightarrow TY \rightarrow \wedge^2L^*
\otimes C^* \rightarrow 0.
$$
\vspace{-.5cm}
\label{suite_tangent_y}
\end{prop}
\rek
The image of $Hom(L,A)_\alpha$ in $T_\alpha Y$ may be described
geometrically, by \cite[theorem 4.3 p.82]{manivel},
as the linear subspace generated by
the tangent vectors to lines through $\alpha$ included in $Y$.\\
\pr
Let $u \in \a$ such that $N(u)=0$.
Let $x=\matttr 100 000 000$ and $y=\matttr 10{\overline u} 000 u00$ be
vectors in $\widehat X$. The tangent spaces $T_x \widehat X$ and 
$T_y \widehat X$ were computed during the proof of proposition
\ref{t_y}; it follows from this computation that
$T_x \widehat X \cap T_y \widehat X = \mattt *.. {R(u)}00 {u^\bot}00$
(the dots replace coefficients above the diagonal which are conjugates
of elements under it)
and that
$\scal{T_x \widehat X , T_y \widehat X} = \mattt *** *0* *{L(u)}*$. 
Note that these spaces do not change when $x$
is replaced by $\lambda .x$, and $y$ by $\nu .y$, $\lambda,\nu \in \C$. 
Therefore, if
$\alpha$ represents the subspace generated by $x$ and $y$, we have
$I_\alpha = T_x \widehat X \cap T_y \widehat X$ and
$S_\alpha = \scal{T_x \widehat X,T_y \widehat X}$. We see that 
$\dim I_\alpha = 3a/2$ and $\dim S_\alpha = 5a/2+2$. The first result
on the ranks of the vector bundles therefore follows.
\lpara

Let $\alpha \in Y$; I now define a map
$T_\alpha Y \rightarrow \wedge^2 L_\alpha^* \otimes C_\alpha^*$.
Let $G(2,V)$ denote the grassmannian of 2-dimensional linear subspaces
of $V$. We use the fact that $T_\alpha Y$, as a subspace of 
$T_\alpha G(2,V)$,
may be described as the set of 
$\varphi : L_\alpha \rightarrow V/L_\alpha$ such
that $\forall x \in L_\alpha,\varphi(x) \in T_x\widehat X/L_\alpha$.
So an element 
$\varphi \in T_\alpha Y \subset Hom(L_\alpha,V/L_\alpha)$ defines a
linear map
$\fonction{\varphi_0}{L_\alpha \otimes L_\alpha}{V^*}
{x \otimes y}{(w \mapsto \det(x,\varphi(y),w)).}$

Now, if $y=\lambda.x$, with 
$\lambda \in \C$, then $\varphi(y) \in T_y \widehat X = T_x \widehat X$, 
and so by lemma \ref{tangent_severi}, $\det(x,\varphi(y),w)=0$ for all
$w \in V$. Therefore, $\varphi_0$ induces a linear map
$\varphi_1:\wedge^2 L_\alpha \rightarrow V^*$.

Moreover, assume there exists $x \in L_\alpha - \{0\}$ such that 
$w \in T_x \widehat X$. Then we have 
$\det(x,\varphi(y),w)=\det(x,w,\varphi(y))=0$ because
$w \in T_x \widehat X$. Choosing  $y \in L_\alpha$ not colinear with
$x$, this proves that
$\varphi_1(\wedge^2 L_\alpha) \subset S_\alpha^\bot$. Since
$S_\alpha^\bot = C_\alpha^*$, we therefore get an element 
$\varphi_2 \in \wedge^2 L_\alpha^* \otimes C_\alpha^*$. The map
$\varphi \mapsto \varphi_2$ is the map 
$T_\alpha Y \rightarrow \wedge^2 L_\alpha^* \otimes C_\alpha^*$ of the
proposition.

From the realization of $T_\alpha Y$ as a subspace of 
$Hom(L_\alpha,V/L_\alpha)$, it is moreover clear that
$Hom(L_\alpha,A_\alpha)$ is a subspace of $T_\alpha Y$.
Assume now that $\varphi_2 = 0$. This implies that if 
$x,y \in L_\alpha$ and $w \in V$, then $\det(x,\varphi(y),w)=0$. By
lemma \ref{tangent_severi} again, this implies that 
$\varphi(y) \in T_x \widehat X$. It follows that 
$\im\ \varphi \subset A_\alpha$ and 
$\varphi \in Hom(L_\alpha,A_\alpha)$.
Since $\dim Y = 25$, the above map 
$\varphi \mapsto \varphi_2$ is surjective and the sequence of the
proposition is exact.
\qed

\lpara

We will see (proposition \ref{wedge2C=A_geometrique})
that the projectivised bundle $\p A$ contains a
subvariety which is isomorphic to the relative grassmannian $G(2,C)$ of
2-dimensional subspaces in $C$. Here is a first result in this
direction.

\begin{prop}
There is a $G$-equivariant injective map of bundles 
$\psi : \wedge^2 C \otimes \wedge^2 L \rightarrow A$. The cokernel bundle is
trivial except when $a=4$, in which case it is a line bundle.
\label{wedge2C=A}
\end{prop}
\pr
Assume first that $a=4$. Let $E$ be a 6-dimensional vector space; we
have already seen that $Y=F(1,3,E)$. So a point $\alpha$ in $Y$
defines a 1-dimensional subspace $E_1$ of $E$ and a 3-dimensional
subspace $E_3$ of $E$; moreover, $E_1 \subset E_3$. Consider now
$E_1,E_3$ as bundles over $Y$.

We have $V \otimes \co_Y= \wedge^2 E$, 
$L = E_1 \wedge E_3 = E_1 \otimes (E_3/E_1)$, 
$A=E_1 \otimes (E/E_3) \oplus \wedge^2 (E_3/E_1)$ and 
$C=\wedge^2(E/E_3)$. Set $A'=E_1 \otimes (E/E_3)$. Recall that if $Z$ is
a 3-dimensional vector space, then $\wedge^2(\wedge^2 Z)$ is
canonically isomorphic with $Z \otimes \wedge^3 Z$. Therefore,
$$
\begin{array}{rcl}
\wedge^2 C \otimes \wedge^2 L & = &
(E/E_3) \otimes \wedge^3(E/E_3) \otimes E_1 \otimes E_1 \otimes
\wedge^2 (E_3/E_1)\\
& = & (E/E_3) \otimes \wedge^3 (E/E_3) \otimes E_1 \otimes \wedge^3 E_3\\
& = & E_1 \otimes (E/E_3) \otimes \wedge^6 E\\
& = & A'.
\end{array}
$$
The last equality follows from the fact that $\wedge^6 E$ is the
trivial line bundle on $Y$. We therefore get the map
$\wedge^2 C \otimes \wedge^2 L \rightarrow A$, which is injective and
has 1-dimensional cokernel.

The case when $a=2$ is similar.

\lpara

Assume now that $a=8$. In this case, I don't know any better proof
than checking the weights. Recall from \cite{bourbaki} the following~:
the highest weight of $V$ is $\lambda = \frac13\poidsesix 456423$ 
and the lowest
is $\frac13\poidsesix {-2}{-4}{-6}{-5}{-4}{-3}$. Let $x$ be a vector of
weight $\lambda$ and $y$ a vector of weight 
$s_{\alpha_1}(\lambda) = \frac13\poidsesix 156423$. We may assume that
$L_\alpha$ is the space generated by $x$ and $y$. 
I claim that the weights
of $C_\alpha$ are
$
\frac13 \poidsesix{-2}{-4}{-6}{-5}{-4}{-3},
\frac13 \poidsesix{-2}{-4}{-6}{-5}{-1}{-3},
\frac13 \poidsesix{-2}{-4}{-6}{-2}{-1}{-3},\\
\frac13 \poidsesix{-2}{-4}{-3}{-2}{-1}{-3}
$ and
$
\frac13 \poidsesix{-2}{-4}{-3}{-2}{-1}{0}.
$
In fact, first, we see that
these weights are obtained from the lowest adding
successively $\alpha_6,\alpha_5,\alpha_4,\alpha_2$ (this proves by the
way that if $L\simeq SL_2 \times SL_5$ is included in 
a Levi factor of $P_3$, then
$C_\alpha$ is an irreducible $SL_5$-module).
Second, the corresponding weight lines are not in $T_x\widehat X$
(resp. neither in $T_y\widehat X$) since the weights of this linear subspace
are the sum of $\lambda$ (resp. $s_{\alpha_1}(\lambda)$) and a
root. Since no root has a coefficient $-3$ in $\alpha_4$, the claim
follows.

Adding the two highest weights of $C_\alpha$ and the two weights of 
$L_\alpha$, one gets
$\frac 13\poidsesix 126423$. This is exactly the highest weight of
$A$. Therefore, there is an $L$-equivariant map 
$\wedge^2 C_\alpha \otimes \wedge^2 L_\alpha \rightarrow
A_\alpha$. Since this is a map between irreducible
$L$-representations, it is also a $P_3$-equivariant map, proving the
proposition.
\qed

\begin{lemm}
Let $\alpha \in Y$ and $x,y \in \p I_\alpha-\p L_\alpha$ such that
$x \equiv y \ \mbox{ mod } \ L_\alpha$. Then $x \in X$ \iff $y \in X$.
\end{lemm}
\pr
Let $z_1 \not = z_2 \in \p L_\alpha$ and $i \in \{1,2\}$.
By definition of $I_\alpha$, 
$x \in T_{z_i}X$. If $x \in X$, then the projective
line $(xz_i)$ through $x$ and $z_i$ 
meets $X$ at the points
$z_i$ and $x$, and with multiplicity at least two at $z_i$. Since $X$
is defined by quadratic equations, $(xz_i) \subset X$. Therefore, the
plane $(xz_1z_2)$ meets $X$ along the three lines 
$(z_1z_2),(xz_1),(xz_2)$; so this plane is included in $X$. Therefore,
$y \in X$.
\qed

\begin{nota}
Let $A' \subset A$ denote the image of $\wedge^2 C \otimes \wedge^2 L$
under the map of proposition \ref{wedge2C=A}.
Let $X(\alpha) \subset \p A'_\alpha$ denote the intersection of the
image of $X$ under the
rational projection $\p I_\alpha \dasharrow \p A_\alpha$ and 
$\p A'_\alpha$.
\end{nota}

\begin{prop}
Assume $a \geq 4$.
Let $\alpha \in Y$ and $x \in X(\alpha)$.
The projectivisation of the inverse of the isomorphism 
$\psi_\alpha : \wedge^2 C_\alpha \otimes \wedge^2 L_\alpha
\rightarrow A'_\alpha$ maps $x$
on the element in $G(2,C_\alpha)$ representing the 2-dimensional space
$T_y \widehat X / S_\alpha \subset V/S_\alpha$, if $y \in I_\alpha$ is
any vector with class $x$ in $\p A_\alpha$.
\label{wedge2C=A_geometrique}
\end{prop}
\pr
Assume first that $a=8$.
Let $\alpha \in Y$. Since $A_\alpha$ is an irreducible
$SL_5$-representation isomorphic with $\wedge^2 \C^5$, there is a
unique non-trivial invariant subvariety in $\p A_\alpha$, and therefore
it is $X(\alpha)$. If $a = 4$, then obviously
we also have $X(\alpha) = \p A'_\alpha$.

Let $a \in \{4,8\}$ and
assume as in the proof of proposition \ref{suite_tangent_y} that
$x=\matttr 100 000 000,y=\matttr 10{\overline u} 000 u00$ and $L_\alpha$
is spanned by $x$ and $y$.
If $z=\matttr 10{\overline v} 000 v00$, with $\scal{u,v} = 0$, then
the 3-dimensional space generated by $x,y,z$ lies in $X$ and 
$T_z \widehat X/S_\alpha \subset C_\alpha$ 
is $2$-dimensional and does not change
if $z$ is replaced by a linear combination of $x,y$ and
$z$. By homogeneity of $X(\alpha)$, 
this fact holds for any $[z] \in X(\alpha)$ and so
we have a well-defined map $X(\alpha) \rightarrow G(2,C_\alpha)$. Since
there is only one such $P_3$-equivariant map, this map also coincides
with the restriction of the projectivisation of $\psi^{-1}$.
\qed

\para

We now assume $a=8$,
and conclude this subsection classifying the $E_6$-orbits in 
$T^*Y$.
By propositions \ref{suite_tangent_y} and \ref{wedge2C=A},
there is a
vector bundle map
$T^* Y \rightarrow Hom(L^* \otimes \wedge^2 L, \wedge^2 C^*)
=Hom(L,\wedge^2 C^*)$; I denote it $h$.

\begin{prop}
Let $\alpha \in Y$ and $f,g \in T^*Y$, 
and assume $f$ and $g$ both
don't vanish. Then $f,g$ lay in the same $E_6$-orbit \iff the two elements
$h(f),h(g) \in Hom(L_\alpha,\wedge^2 C_\alpha^*)$
lay in the same $(GL(L_\alpha) \times GL(C_\alpha))$-orbit.
\label{orbite_tangent}
\end{prop}
\noindent
In view of lemma \ref{orbites}, this gives a complete understanding of
the $E_6$-orbits in $T^*Y$.\\
\pr
Let $P \subset E_6$ be the stabilizor of $\alpha$ and $L(P)$ a Levi
factor of $P$. We know that the image of $L(P)$ in
$End(Hom(L_\alpha, \wedge^2 C_\alpha))$ is the same as that of
$GL(L_\alpha) \times GL(C_\alpha)$.
If $h(f)=h(g)=0$, then, by proposition \ref{suite_tangent_y},
$f$ and $g$ are elements in
$(\wedge^2 L_\alpha \otimes C_\alpha)-\{0\}$,
which is obviously homogeneous under
$L(P)$, and so lay in the same $P$-orbit.

Assume $h(f) \not = 0$ and $h(g) \not = 0$. Since by hypothesis $h(f)$
and $h(g)$ lay in the same $L(P)$-orbit, we may
assume that $h(f)=h(g)$. Let $R_u(P)$ denote the unipotent radical of
$P$; $R_u(P)$ acts trivially on the irreducible $P$-representation
$Hom(L_\alpha, \wedge^2 C_\alpha^*)$.
Therefore, it is enough to prove that the $R_u(P)$-orbit
of $f$ is dense in $h^{-1}(h(f))$. Equivalently, we will prove that 
the image of the action of the Lie algebra of $R_u(P)$ on $f$ contains
$\wedge^2 L_\alpha \otimes C_\alpha$.

It is enough to prove this when $h(f)$ is in the minimal non-zero
orbit of $GL(L_\alpha) \times GL(C_\alpha)$ in
$Hom(L_\alpha, \wedge^2 C_\alpha^*)$.
This, in turn, can be verified at the level of
weights. In fact, we assume that $h(f)$ is a highest weight vector of
$Hom(L_\alpha, \wedge^2 C_\alpha^*)$. Therefore, $h(f)$ has weight
$\poidsesix 112212$. In fact, as we saw in the proof of proposition
\ref{wedge2C=A}, the highest weight
of $L^*_\alpha$ is
$\frac13 \poidsesix{-1}{-5}{-6}{-4}{-2}{-3}$, and
the two highest weights of $C^*_\alpha$ are
$\frac13 \poidsesix{2}{4}{6}{5}{4}{3}$ and
$\frac13 \poidsesix246513$.

Since the weight of $\wedge^2L$ is
$\frac13 \poidsesix{5\ }{10}{\ 12\ }{8\ }{4}{6}$ and the heighest
weight of $C$ is $\frac13 \poidsesix{-2\ }{-4}{\ -3\ }{-2\ }{-1}{0}$,
the heighest weight of $\wedge^2 L_\alpha \otimes C$ is
$\poidsesix 123211 = \omega_2$.
Since this is the
highest weight of
$Hom(L_\alpha, \wedge^2 C_\alpha^*)$
plus $\alpha_3 + \alpha_4$, which is a root of $R_u(P)$, we are done.
\qed

%******************************************************************************

%******************************************************************************

%******************************************************************************

\sectionplus{Mukai flops of type $E_6$}

\label{flop}

Let $(V,X)$ be a Scorza variety and $(V^*,X^\dual)$ the dual Scorza
variety.
An element in $T^*X$
will be denoted $(x,\alpha)$, where $x \in X$ and $\alpha$ is a
linear form on $T_xX$.
The flop is a map $T^*X \dasharrow X^\dual,(x,\alpha) \mapsto 
(h,\kappa)=(h(x,\alpha),\kappa(x,\alpha))$.

For flops of type $E_{6,I}$, the element $h(x,\alpha)$ 
was described in the preceeding
section. The complete description of Mukai flops should also include a
formula for $\kappa(x,\alpha)$. However, it is not easy to follow the
identification of $T^*_xE_6/P_1$, seen more or less as a subspace of
$V^*$, with a subspace of the Lie algebra of
$G$. Instead, given $h(x,\alpha)$, I explain in the next 
subsection a general geometric way to put
our hands on $\kappa(x,\alpha)$.

\subsectionplus{Canonical isomorphism of quotients of tangent spaces
to flag varieties} 

\label{general}

Let $G$ be a reductive algebraic group and let
$\cal P,Q$ denote two flag varieties parametrizing two classes of
parabolic subgroups of $G$. Let $\cal R$ denote the flag variety of
parabolic subgroups which are intersections of a parabolic subgroup in
$\cal P$ and a parabolic subgroup in $\cal Q$.
Since a parabolic subgroup in $\cal R$ is contained in
exactly one subgroup in $\cal P$ (resp. $\cal Q$),
$\cal R$ is canonically isomorphic
with a subvariety of ${\cal P} \times {\cal Q}$; an element in 
$\cal R$ will threrefore be denoted $(x,y)$, with $x \in \cal P$ and
$y \in \cal Q$.

If $x \in \cal P$, let ${\cal Q}_x$ denote the variety of parabolic
subgroups $y$ such that $(x,y) \in \cal R$, and define similarly
${\cal P}_y$. The following quite general theorem allows, as a special
case, describing Mukai flops $T^*X \dasharrow T^*X^\dual$
of type $E_{6,I}$
as soon as we know the composition
$T^*X \dasharrow T^*X^\dual \rightarrow X^\dual$. Since it is the case
by theorem \ref{cotangent}, it will be easy to deduce a Jordan-theoretic
formula for this Mukai flop, in proposition \ref{flop_premier}.

\lpara

In the next theorem, $(C,0,t)$ will denote a pointed curve $(C,0)$
which is smooth at $0$,
together with a tangent vector $t$ at the point $0$. Moreover, if 
$Y$ is an algebraic variety and $f:C \rightarrow Y$ is a map, 
then $f'(0) \in T_{f(0)}Y$ 
will denote the derivative $df_0(t)$.
I say that there is a Mukai flop
$T^*{\cal P} \dasharrow T^*{\cal Q}$
if the natural maps
$T^*{\cal P} \rightarrow \g,T^*{\cal Q} \rightarrow \g$ are birational
and have the same image.

\begin{theo}
Let $(x,y) \in {\cal R}$. Then there is a canonical isomorphism
$\mu(x,y):\frac{T_x{\cal P}}{T_x{\cal P}_y}
\rightarrow \frac{T_y{\cal Q}}{T_y{\cal Q}_x}$. If
$(C,0,t)$ is as above, and if
$\gamma:(C,0) \rightarrow ({\cal P},x)$ is any map,
such an
isomorphism maps the class of $\gamma'(0)$ on the class of 
$\delta'(0)$, if $\delta:(C,0) \rightarrow ({\cal Q},y)$ is any map
such that $(\gamma,\delta)(C) \subset {\cal R}$.\\
If, moreover, there is a Mukai flop 
$T^*{\cal P} \dasharrow T^*{\cal Q}$, then this flop maps a generic form
$f \in (T_x{\cal P}/T_x{\cal P}_y)^*$ to
$(y,\tr \mu(x,y)^{-1}(f))$.
\label{mukai}
\end{theo}
\pr
Let $\pi_P:{\cal R}\rightarrow {\cal P}$ and
$\pi_Q:{\cal R}\rightarrow {\cal Q}$ denote the natural projections.
Consider the diagram
$$
\begin{array}{ccccc}
\frac{T_x{\cal P}}{T_x{\cal P}_y} & \stackrel{\varphi_P}{\leftarrow}  &
\frac{T_{(x,y)}{\cal R}}{\scal{T_{(x,y)}\pi_P^{-1}(x),T_{(x,y)}\pi_Q^{-1}(y)}}
& \stackrel{\varphi_Q}{\rightarrow} & \frac{T_y{\cal Q}}{T_y{\cal Q}_x}\\
\parallel & & \parallel & & \parallel \vspace{.15cm}\\
\frac{\g/\plie}{\q/\plie} & \simeq &
\frac{\g/(\plie \cap \q)}{\scal{\plie,\q}/(\plie \cap \q)} & \simeq &
\frac{\g/\q}{\plie/\q}\ \ ,
\end{array}
$$
where $\varphi_P$
(resp. $\varphi_Q$) is induced by the
differential $d_{(x,y)} \pi_P$ (resp. $d_{(x,y)}\pi_Q$).
All the terms on the second line are canonically isomorphic
with $\g/\scal{\plie,\q}$.
Obviously, the diagram commutes, so $\varphi_P$ and $\varphi_Q$ are
isomorphisms.
Let $\mu(x,y):\frac{T_x{\cal P}}{T_x{\cal P}_y}
\rightarrow \frac{T_y{\cal Q}}{T_y{\cal Q}_x}$ denote the canonical
isomorphism $\varphi_Q \circ \varphi_P^{-1}$.

Let $f \in (T_x{\cal P}/T_x{\cal P}_y)^*$ be generic.
Let $(y',f') \in T^*{\cal Q}$ denote the image of $(x,f)$ by the flop
$T^*{\cal P} \dasharrow T^*{\cal Q}$. By \cite{dual}, $y'$ is the only
element in $\cal Q$ such that $f$ vanishes on $T_x{\cal P}_y$;
since by assumption, $f$ vanishes on $T_x{\cal P}_y$, $y'=y$.
Moreover, we have canonical isomorphisms
$(T_x{\cal P}/T_x{\cal P}_y)^* \simeq \u(\plie) \cap \u(\q) \simeq
(T_y{\cal Q}/T_y{\cal Q}_x)^*$, and under this isomorphism, $f$ is
mapped to $f'$ by definition of the Mukai flop. It is clear that this
isomorphism is the transpose of $\mu(x,y)^{-1}$, so the last claim of
the proposition is proved.

\lpara

If $(\gamma,\delta)$ are as in the proposition, then 
$(\gamma'(0),\delta'(0)) \in T_{(x,y)} {\cal R}$; denote by
$[\gamma'(0),\delta'(0)]$ its class in
$\frac{T_{(x,y)}{\cal R}}
{\scal{T_{(x,y)}\pi_P^{-1}(x),T_{(x,y)}\pi_Q^{-1}(y)}}$. By definition
of $\varphi_P$ and $\varphi_Q$, we have
$\varphi_P([\gamma'(0),\delta'(0)]) = [\gamma'(0)]$ and
$\varphi_Q([\gamma'(0),\delta'(0)]) = [\delta'(0)]$. We therefore
have, as expected,
$\varphi_Q \circ \varphi_P^{-1}([\gamma'(0)] = [\delta'(0)] $.

\qed

\subsectionplus{Mukai flop of type $E_{6,I}$ in terms of Jordan algebras}

Let $\mu(x,y)$ denote the isomorphism of proposition \ref{mukai}.
In this subsection, I give an expression of $\mu(x,y)$
in the case of Scorza varieties, in terms of Jordan
algebras. Therefore, this gives also a formula for the Mukai flop.

More precisely, let $(V,X)$ be a Scorza variety of type $(n,a)$ and let
$(V^*,X^\dual)$ be the dual Scorza variety.
Let $(x,h) \in X \times X^\dual$ such that $x \vdash h$.
Let us choose $(\tilde x, \tilde h) \in V \times V^*$ such that 
$[\tilde x] = x$ and $[\tilde h] = h$. This identifies $T_xX$
(resp. $T_hX^\dual$) with $T_{\tilde x} \widehat X/\C.\tilde x$
(resp. $T_{\tilde h} \widehat {X^\dual}/\C.\tilde h$). The previous
isomorphism $\mu(x,h) : T_xX / T_xC_h \simeq T_hX^\dual / T_hC_x$ induces an
isomorphism $T_{\tilde x} \widehat X/T_{\tilde x} \widehat{C_h}
\simeq T_{\tilde h} \widehat {X^\dual} / T_{\tilde h} \widehat {C_x}$.

The goal of this section is to give a formula for this isomorphism in
Jordan terms.

\lpara

For
$A,B \in V$, let
$\sigma_A(B) \in V^*$ denote the linear form $U \mapsto D^2_A \det(B,U)$.
Note that this is equal, modulo $D_A \det$, to $S_A(B)$ \cite{dual}.
For $h \in X^\dual$, let $V(h):={(T_hX^\dual)}^\bot \subset V$.

\begin{lemm}
Let $A \in V(h)$. Then $\sigma_A(\tilde x)$ is proportional to $\tilde h$.
\label{proportionnel}
\end{lemm}
\pr
We can assume that $V=H_3(\a)$ and $X \subset \p V$ is the variety of
rank one elements. Identify $V$ and $V^*$ as usually.
Since $X^\dual$ is homogeneous under $G$, we may assume that
$h=\matttr 000 000 001$. Then 
$V(h) = \left \{ \matttr **0 **0 000 \right \}$.
It is enough to prove the lemma for generic $A$ in $V(h)$, so we
may assume that $A$ has rank $2$. Moreover, let $G_h$ denote the stabilizor
in $Aut(X)$ of $h$. It is clear that $G_h$ acts transitively on the
set of rank $2$ elements of $V(h)$, and on the set of its 
rank 1 elements. So we may assume $\tilde x=\matttr100 000 000$. 
Moreover, for the stabilizor
of $x$ in $G_h$, the set of $A$'s of rank $2$ is made of two orbits and
$A = \matttr 100 010 000$ is in the open orbit, as it is easily checked
case by case.

Therefore, it is enough to compute $\sigma_A(\tilde x)$ for these choices
of $x$ and $A$. Let $m=(m_{i,j}) \in V$. Then one computes
\begin{equation}
D^2_A \det(m,m) = \sum_{i<n} m_{i,i}m_{n,n} - \sum_{i<n} N(m_{n,i}).
\label{derivee_seconde}
\end{equation}

The lemma immediately follows.
\qed

\lpara
Let as before $x \in X,h\in C_x,A \in V(h)$, and let 
$\tilde x \in V,\tilde h \in V^*$ represent $x$ and $h$.  

\begin{prop}
Let $v \in T_{\tilde x} \widehat X$, and let $[v]$ denote its class
in $T_{\tilde x} \widehat X/T_{\tilde x} \widehat{C_h}
\simeq T_xX / T_xC_h$.
If $\sigma_A(\tilde x) = \tilde h$, then the vector
$$[\sigma_A(v)] \in T_{\tilde h}\widehat{X^\dual} / T_{\tilde h}
\widehat{C_x} \simeq T_hX^\dual / T_hC_x$$ 
identifies with $\mu(x,h)([v])$.
\label{flop_premier}
\end{prop}
\noindent
The isomorphism $T_{\tilde x} \widehat X/T_{\tilde x} \widehat{C_h}
\simeq T_xX / T_xC_h$ depends on the choice of $\tilde x$, and the
isomorphism $T_{\tilde h}
\widehat{C_x} \simeq T_hX^\dual / T_hC_x$ depends on the choice of
$\tilde h$. However, the proposition says that the corresponding map
$T_xX / T_xC_h \rightarrow T_hX^\dual / T_hC_x$ does not depend on these
choices, neither on the choice of $A$, as long as
$\sigma_A(\tilde x) = \tilde h$.

\lpara
\noindent
\pr
As in the previous lemma, we assume that $V = H_n(\a)$.
Let $X_r$ denote the variety of rank $r$ matrices. If $B \in X_{n-1}$,
then $D_B \det$ belongs to $X^\dual$. Since $D_A \det  = h$ and
$\widehat{T_xX} \subset \widehat{T_AX_{n-1}}$, we have the implication
$u \in \widehat{T_xX} \Longrightarrow \sigma_A(u) \in
\widehat{T_hX^\dual}$.
\lpara

Now, let $v \in \widehat{T_xX}$ and let $u$ be the class of $v$ in
$\widehat{T_xX}/\C.\tilde x$. Let $\varphi(\tilde x,\tilde h,A)(u)$
denote the element of $T_hX^\dual$ corresponding to the class of
$\sigma_A(v)$ in $\widehat{T_hX^\dual}/\C.\tilde h$ (by lemma
\ref{proportionnel}, this
class depends only on $u$).

We first show that if $\tilde x,\tilde h,A$ are multiplied by a
scalor, then $\varphi(\tilde x,\tilde h,A)$ does not vary. So let
$\lambda,\mu,\nu \in \C$, and assume 
$\sigma_{\nu.A}(\lambda.\tilde x) = \mu.\tilde h$. Since
by assumption $\sigma_A(\tilde x) = \tilde h$, this means that
$\nu^{n-2}.\lambda = \mu$.
Now, $\lambda.\tilde x$ will identify $u$ with $\lambda.v$. Then,
$\sigma_{\nu.A}(\lambda.v) = \nu^{n-2}.\lambda.\sigma_A(v) 
= \mu.\sigma_A(v)$, and
the class of this vector will identify with 
$\varphi(\tilde x,\tilde h,A)(u) \in T_hX^\dual$ with the choice 
$\mu.\tilde h$ instead of $\tilde h$.

\lpara

Therefore, the claim is proved, and one can choose the same elements 
$\tilde h,\tilde x,A$ as in the proof of the lemma.
By formula (\ref{derivee_seconde}), if 
$v=\matttr t{\overline w}{\overline z}
w00
z00 \in \widehat{T_xX}$,
then $\sigma_A(v)=\matttr 00{-\overline z}
000
{-z}0t$, if one
identifies $V$ and $V^*$ via the usual scalar product.
Note that $T_{\tilde h}\widehat C_x = \left \{
\matttr 000 00* 0** \right \}$, so that the class of $\sigma_A(v)$ does not
depend on $w$, neither on $t$.
\lpara

Let us now compute $\mu(x,h)([v])$ using
theorem \ref{mukai}, and assuming $w=0$.
So let $z \in \a$ and $t \in \C$.
Recall that generic elements of $\widehat X$
can be written as $\nu_2(\alpha,z_1,z_2)$, with 
$\alpha \in \C,z_1,z_2 \in \a$, and 
$\nu_2(\alpha,z_1,z_2) = \matttr {\alpha^2}{\alpha\overline z_1}
{\alpha\overline z_2}{\alpha z_1}{N(z_1)}{z_1\overline z_2}
{\alpha z_2}{z_2 \overline z_1}{N(z_2)}$.
Denote $x(t)=\nu_2(1,0,tz)= \matttr 1 0 {\overline zt} 0 0 0 {zt} 0
{N(z)t^2}$; we
have $x(t) \in \widehat X$ and 
$x'(0) = \matttr 0 0 {\overline z} 0 0 0 z 0 0$.

Differetiating $\nu_2$ we get
$$\widehat{T_{x(t)}X} = 
\left \{
\mattt {2\alpha}{\overline z_1}{\overline z_2 + \alpha t \overline z}
{z_1}{0}{t z_1 \overline z}
{\alpha t z + z_2}{tz\overline z_1}{t(z\overline z_2 + z_2 \overline z)}
: \alpha \in \C,z_1,z_2 \in \a \right \}
$$

Recall that the incidence relation between $X$ and $X^\dual$ is :
$x \vdash h$ if $h \supset \widehat{T_xX}$; therefore, if we set
$h(t)=\matttr{t^2N(z)/2} 0 {-t\overline z} 0 0 0 {-tz} 0 1$, we have
$x(t) \vdash h(t)$, and since
$h'(0)=-\matttr 00{\overline z} 000 {z}00$, the proposition follows.
\qed

\subsectionplus{Mukai flops for Scorza varieties in terms of $\a$-blow-up}

\label{eclate}

The simplest Mukai flop $T^* \p^n \dasharrow T^* {(\p^n)}^\dual$ can be
resolved blowing-up the zero section. Let's recall this
construction. Let $Z \subset T^* \p^n$ be the zero section, and let
$B$ be the blow-up of $T^* \p^n$ along $Z$. It is known that there is
a map $B \rightarrow T^* {(\p^n)}^\dual$ such that the following triangle
commutes~:
\begin{equation}
\label{commutatif}
\begin{array}{ccccc}
 && B \\
& \swarrow & & \searrow\\
T^* \p^n & & \dasharrow & & T^* {(\p^n)}^\dual.
\end{array}
\end{equation}

Moreover, this is the minimal resolution, in the sense that for any other
$B'$ with the same property, there is a map $B' \rightarrow B$ and an
obvious commutative diagram.

\lpara

In this subsection, I give a similar resolution of the rational map
$T^*X \dasharrow T^*X^\dual$, if $X$ is a Scorza variety and $X^\dual$
the corresponding dual Scorza variety. In fact, the main idea is that
since $X$ behaves like a projective space $\p^n_\a$ over $\a$, one
should consider an ``$\a$-blow-up''.

Let me make a heuristic comment.
Given a composition algebra $\a$, I believe in the existence of
a category $\a-Var$ of $\a$-varieties, containing projective
spaces and grassmanians over $\a$. Moreover, if $Y \subset X$ is
a closed immersion in this category,
then there should be an object $Bl_X(Y)$ over $X$
defined by a universal property analogous to that defining usual
blow-ups, but in the category $\a-Var$.
Since for the moment I don't know how to define $\a-Var$, I will not
give this construction here. In the following we
will only have very simple $\a$-blow-ups to do, and in these
simple cases we can guess what the blow-up should be.

\lpara

So let $\a$ be a composition algebra over $\C$ of dimension $a$ 
and $n \geq 2$ an integer, with
$n=2$ if $\a = \oc$. Let the affine space $\A^n_\a$ be 
just $\a^n$, the affine
$(an)$-dimensional space over $\C$.

Recall that in subsection \ref{p1o}, I introduced a map 
$\overline \nu_2:\A^2_\O \dasharrow \p^1_\O$, where by definition
$\p^1_\O$ is an 8-dimensional smooth quadric.
Recall also the rank 8 vector bundle $L$ over $\p^1_\O$ of
proposition \ref{ls}. By definition, $L$ is a
subbundle of the trivial bundle $\A^2_\O \otimes \co_{\p^1_\O}$ of rank
16 over $\p^1_\O$. Therefore, if $\cal L$ denotes the total space of
the vector bundle $L$, there is an inclusion
${\cal L} \subset \A^2_\O \times \p^1_\O$. Therefore we have a map
${\cal L} \rightarrow \A^2_\O$.

The case of associative algebras is simpler and was studied in
\cite{projective}~: recall that there is a rational map
$\overline \nu_2 : \A^n_\a \dasharrow \p^{n-1}_\a$ and a rank $a$ subbundle
$\cal L$ of the trivial bundle $\A^n_\a \otimes \co_{\p^{n-1}_\a}$ with
fiber $\A^n_\a$ over $\p^{n-1}_\a$. This subbundle is also defined by
${\cal L}_x = \overline{ \{ v \in \A^n_\a : 
\overline \nu_2(v) \mbox{ is defined and } \overline \nu_2(v) = x \} }$.
Let $\cal L$ denote its total space.

\begin{defi}
The $\a$-blow-up $Bl_{\A^n_\a}(0)$ of the affine space
$\A^n_\a$ at the origin is the map
${\cal L} \rightarrow \A^n_\a$.
\end{defi}

Recall that in $\A^2_\O$ there are three $Spin_{10}$-orbits~: the open
orbit, the point $0$, and the affine cone $\widehat \S$ over a spinor
variety $\S \subset \p \A^2_\O$.
The map $Bl_{\A^2_\O}(0) \rightarrow \A^2_\O$ is an isomorphism above
the open orbit, and the fiber over $0$ is isomorphic with
$\p^1_\O$. Except for the existence of the intermediate 
orbit $\widehat \S-\{0\}$ in
$\A^2_\O$, the situation is therefore very similar to that of the
usual blow-up of the origin in $\A^2_\C$. A similar statement holds
in general for the blow-up of $\A^n_\a$. The following result gives
another analogy with usual blow-ups~:

\begin{prop}
This $\a$-blow-up is the minimal resolution of the rational map
$\overline \nu_2 : \A^n_\a \dasharrow \p^{n-1}_\a$.
\label{resolution}
\end{prop}
\pr
Let $\pi : {\cal L} \rightarrow \A^n_\a$ denote this $\a$-blow-up.
The restriction of $\pi$ to the regular locus of $\overline \nu_2$ is
an isomorphism.
By definition, there are maps 
$Bl_{\A^n_\a}(0) \rightarrow \A^n_\a$ and
$Bl_{\A^n_\a}(0) \rightarrow \p^{n-1}_\a$ such that the diagram
$$
\begin{array}{ccccc}
 && Bl_{\A^n_\a}(0) \\
& \swarrow & & \searrow\\
\overline \nu_2 : \A^n_\a & & \dasharrow & & \p^{n-1}_\a
\end{array}
$$
commutes.

Let
$$
\begin{array}{ccccc}
 && B' \\
& \swarrow & & \searrow\\
\overline \nu_2 : \A^n_\a & & \dasharrow & & \p^{n-1}_\a
\end{array}
$$
be another resolution. Then we have a map 
$B' \rightarrow \A^n_\a \times \p^{n-1}_\a$. Since the above diagram is
commutative, the image of this map is 
$Bl_{\A^n_\a}(0)$, and we get the desired map
$B' \rightarrow Bl_{\A^n_\a}(0)$.
\qed

\Para

Let $(V,X)$ be a Scorza variety of type $(n,a)$ with $n \geq 3$.
The above construction of the blow-up of a point in the fixed
vector-space $\A^{n-1}_\a$ extends readily to a blow-up of the zero
section in the vector bundle $T^*X$. In fact, let $x \in X$;
recall (theorem \ref{cotangent}) that we have a
quadratic map $\nu_x^- : T_x^*X \rightarrow N_x^*X$, where 
$N_xX = V/\widehat{T_xX}$, and so a rational
map $\overline {\nu_x^-} : T_x^*X \dasharrow \p  N_x^*X$.
This map is isomorphic with our model map
$\overline \nu_2 : \a^{n-1} \dasharrow \p^{n-2}_\a$.
Letting $x$
vary, we get an algebraic map
$\nu^- : T^*X \rightarrow N^*X$ over $X$, and so a rational map
$\overline {\nu^-} : T^*X \dasharrow \p N^*X$. 

In subsection \ref{variety}, the projectivisation of
the image of $\nu^-$ was denoted $\p N(X)^\dual$;
$\p N(X)^\dual$ is a locally trivial
fibration over $X$ with fibers Scorza varieties of type $(n-1,a)$.
Let $p_X : \p N(X)^\dual \rightarrow X$ denote the natural projection. 

Consider the bundle $p_X^* T^*X$ above $\p N(X)^\dual$. An element of this
bundle will be denoted $(x,h,f)$, with $x \in X,h \in \p N(X)_x^\dual$ and 
$f \in T_x^*X$.
Globalizing the above construction, let
${\cal L} \subset p_X^* T^*X$ be defined as the closure of the set of
$(x,h,f) \in p_X^* T^*X$ such that $\overline{\nu_x^-}(f)$
is defined and equals $h$.

\begin{lemm}
${\cal L} \subset p_X^* T^*X$ is a subbundle.
\end{lemm}
\pr
Assume first that $a=8$.
Then it is simply a global version of proposition \ref{ls}. 
By theorem \ref{cotangent},
We know
that $\nu^-$
is a global algebraic map, which on each fiber $T^*X$ is isomorphic
with the map $\nu_2^- : \ok \oplus \ok \dasharrow \p^1_\O$ defined in
subsection \ref{p1o}.
Therefore, the argument of proposition \ref{ls} works in this
situation.
The case of associative composition algebras $\a$ is similar and left
to the reader.
\qed

Let $Z \subset T^*X$ denote the zero section.

\begin{defi}
The $\a$-blow-up $Bl_{T^*X}(Z)$ of $T^*X$ along $Z$ is the map
${\cal L} \rightarrow T^*X$.
\end{defi}

\begin{theo}
This $\a$-blow-up is the minimal resolution of the Mukai flop
$\mu : T^*X \dasharrow T^*X^\dual$.
\label{eclatement}
\end{theo}
\pr
Globalizing the proof of proposition \ref{resolution}, we see that
$Bl_{T^*X}(Z)$ is the minimal resolution of the rational map 
$\overline {\nu^-} : T^*X \dasharrow \p N(X)^\dual$.
In view of theorem
\ref{cotangent}, it is also the minimal resolution of the composition
$T^*X \stackrel{\mu}{\dasharrow} T^*X^\dual \rightarrow X^\dual$.
Now, by theorem \ref{mukai}, resolving the Mukai flop
$T^*X \dasharrow T^*X^\dual$ is equivalent with resolving its
projection to $X^\dual$, so the theorem follows.
\qed

\subsectionplus{Mukai flop of type $E_{6,II}$}

Let $Y = E_6 / P_3$ be the homogeneous space considered in
subsection \ref{para_tangent_severi}, $Y^\dual = E_6 / P_5$ the ``dual''
homogeneous space and $A,B,C$ the homogeneous vector bundles over $Y$ defined
there. Let also $X = E_6 / P_1$ and $X^\dual = E_6 / P_6$.

We already used the fact that $Y$ is isomorphic with the Fano variety
of projective lines included in $X$. Similarly, $Y^\dual$ is the
variety of lines included in $X^\dual$. Denote as before $\p V$ the
ambient space of $X$. As we have already
seen, $X^\dual$ identifies with the set of hyperplanes in $V$ which
contain two tangent spaces to $X$.

Therefore, given a point $\alpha \in Y$, which represents a projective
line
$l_\alpha$ contained in $X$, and given two points $x\not = y \in X$,
any hyperplane $h \subset \p V$ which contains the span of $T_xX$ and
$T_yX$ can be considered as an element of $X^\dual$. A codimension two
subspace $V_\beta \subset V$ containing this span defines a pencil of 
hyperplanes belonging to $X^\dual$, or a point in $Y^\dual$.

Let $\alpha \in Y$. Recall that the linear space
$V/\scal{\widehat{T_xX},\widehat{T_yX}}$ (where $x$ and $y$ are
different points of the line $l_\alpha$) was denoted $C_\alpha$ in
subsection \ref{para_tangent_severi}. Let $\p C^*$ denote 
the projective bundle
over $Y$ and $G(3,C)$ the relative
grassmannian of 3-spaces in $C$. The preceeding remarks show that
there are natural maps $\p C^* \rightarrow X^\dual$ and
$G(3,C) \rightarrow Y^\dual$.
Let $f:G(3,C) \rightarrow Y^\dual$ be this map. 

\lpara

For any
$\alpha \in Y$, let 
$g_\alpha : Hom(L_\alpha,\wedge^2 C_\alpha^*)
\dasharrow G(3,C_\alpha)$ be the map defined by lemmas \ref{conditions}
and \ref{existence}
using $F = C_\alpha$ (namely, 
$g_\alpha(\varphi) = U(\varphi(l_1),\varphi(l_2))$ for 
$\varphi \in Hom(L_\alpha,\wedge^2 C_\alpha^*)$ 
and any non-colinear $l_1,l_2 \in L_\alpha$).
By propositions \ref{suite_tangent_y} and \ref{wedge2C=A}, there is a natural
vector bundle map
$T^* Y \rightarrow Hom(L^* \otimes \wedge^2 L,
\wedge^2 C^*) = Hom(L,\wedge^2 C^*)$, which I denote $h$.

Let finally $\mu : T^* Y \dasharrow T^* Y^\dual$ be the Mukai flop and 
$\pi : T^* Y^\dual \rightarrow Y^\dual$ the structure map.

\begin{theo}
The composition
$$
T^*Y \stackrel{h}{\rightarrow}
Hom(L,\wedge^2 C^*)
\stackrel{g}{\dasharrow}
G(3,C) \stackrel{f}{\rightarrow}
Y^\dual
$$
equals the composition
$$
T^*Y \stackrel{\mu}{\dasharrow}
T^*Y^\dual \stackrel{\pi}{\rightarrow}
Y^\dual.
$$
\label{flop_second}
\end{theo}
\rek
This describes the rational map $\pi \circ \mu$. The rational
map $\mu$ itself is then
described using proposition \ref{mukai}.\\
\pr
Let $\alpha \in Y$ and generic $\eta \in T_\alpha^* Y$.
We know that $\pi \circ \mu (\eta)$ is the unique 
$\beta \in f(G(3,C_\alpha))$ such that $\eta$ vanishes on the tangent
space $T_\alpha SC_\beta$ at $\alpha$ of the Schubert cell 
$SC_\beta \subset Y$
defined by $\beta$. So first, we compute $T_\alpha SC_\beta$.

\lpara

If $\beta = f(\beta_0)$, with $\beta_0 \in G(3,C_\alpha)$, let
$c_\beta \subset C_\alpha$ be the 3-dimensional subspace
corresponding to $\beta_0$. Let $a_\beta$ denote the image of
$\wedge^2 c_\beta \otimes \wedge^2 L_\alpha
\subset \wedge^2 C_\alpha \otimes \wedge^2 L_\alpha$ in
$A_\alpha$ under the isomorphism of proposition
\ref{wedge2C=A}.
By proposition \ref{wedge2C=A_geometrique},
if the class modulo $L_\alpha$ of $x \in X$ is in 
$\p a_\beta$,
then $\widehat{T_xX}/S_\alpha \subset c_\beta$. 

Let
$v_\beta \subset V$ denote the inverse image of $a_\beta$ under the
projection $I_\alpha \rightarrow A_\alpha = I_\alpha/L_\alpha$.
Since
$SC_\beta$ is the variety of lines $l \subset X$ such that
$\forall x \in l,\widehat{T_xX}/S_\alpha \subset c_\beta$, we deduce
$G(2,v_\beta) \subset SC_\beta$.

Now, given $\alpha \in Y$, the cell $SC_\alpha \in Y^\dual$ identifies
with $G(3,C_\alpha)$.
By symmetry, $SC_\beta$ is also isomorphic
with a six-dimensional grassmannian, and so $G(2,v_\beta) = SC_\beta$.
Therefore, $T_\alpha SC_\beta = Hom(L_\alpha,a_\beta)$.

\lpara

Now, we complete the proof. We already saw that a
cotangent form $\eta \in T^*_\alpha Y$ defines an element
$h(\eta) \in Hom(L_\alpha,\wedge^2 C_\alpha^*)$. Given the previous
computation of $T_\alpha SC_\beta$, $\eta$ vanishes on $T_\alpha SC_\beta$
\iff $\wedge^2 c_\beta \bot \im\ h(\eta)$. Therefore, we can conclude
thanks to lemmas \ref{conditions} and \ref{existence}.
\qed

\lpara

Recall that if $V_2$ and $V_5$ are vector spaces of respective dimensions
2 and 5, there are 8 $(GL(V_2) \times GL(V_5))$-orbits in
$Hom(V_2,\wedge^2V_5)$, which were given a
label in lemma \ref{orbites}. We also use the standart labels for
nilpotent orbits, as in \cite[p.202]{mcgovern}.

\begin{coro}
Let $\alpha \in Y$ and $0 \not = f \in T^*_\alpha Y$.
Under the natural map
$T^*Y \rightarrow \mathfrak{e}_6$, $f$ is mapped to the nilpotent orbit 
with the same label as that of the $(GL(L) \times GL(C))$-orbit
$h(f) \in Hom(L, \wedge^2 C^*)$ belongs. The Mukai
flop is defined exactly on the open orbit of $T^*Y$.
\label{lieu}
\end{coro}
\noindent
In this corollary, I mean that if $h(f)$ is the orbit labelled 
$3A_{1,a},3A_{1,b}$ or $3A_{1,c}$, then it is mapped on the nilpotent orbit
labelled $3A_1$.\\
\pr
This corollary follows from dimension arguments, which are not so
illuminating on the geometry of the resolution.
The map $T^*Y \rightarrow \mathfrak{e}_6$ being birational and proper, it has
50-dimensional closed image; so it is the closure of the unique
50-dimensional orbit in $\mathfrak e_6$, labelled $A_2 + 2A_1$.
By the given graph of orbit closures \cite[p.212]{mcgovern}, the
image of $T^*Y$ is the union of the orbits labelled
$0,A_1,2A_1,3A_1,A_2,A_2+A_1,A_2+2A_1$.

Let $\alpha \in Y$ and $f,g \in T^*_\alpha Y$, 
with $f \not = 0$ and $g \not = 0$. 
Then $f$ and $g$ lay in the same
$E_6$-orbit \iff 
$h(f),h(g) \in Hom(L_\alpha, 
\wedge^2 C_\alpha^*)$ lay in
the same $(GL(L_\alpha) \times GL(C_\alpha))$-orbit, by proposition
\ref{orbite_tangent}. It is clear that the zero section in $T^*Y$ is
mapped to the 0-orbit in $\mathfrak e_6$. 
Let us label the other $E_6$-orbits in
$T^*Y$ by the labels of their images in $Hom(L,\wedge^2C^*)$.

We first begin with a trivial remark~: the image of an orbit in
$T^*Y$ is an orbit in $\mathfrak e_6$ of non-greater dimension. From
this it follows that the orbits in $T^*Y$ labelled 
$A_2+2A_1,A_2+A_1,A_1$ map to the orbits in $\mathfrak e_6$ with the same
label.

Suppose the orbit in $T^*Y$ labelled $3A_{1,a}$ maps to the orbit in
$\mathfrak e_6$ labelled $A_2$. The the fibers above the nilpotent
orbit $A_2$ would have dimension 3, and the preimage of the nilpotent
orbit labelled $3A_1$ would be included in the orbits labelled
$3A_{1,b}$ and $3A_{1,c}$. So the fibers above this orbit would have
dimension 1 or 2, contradicting the semi-continuity of the dimensions
of the fibers of a morphism. Therefore, $3A_{1,a}$ maps to $3A_1$.

We know that the resolution $T^*Y \rightarrow \mathfrak e_6$ is
semi-small, so the 38-dimensional orbit in $T^*Y$ labelled $2A_1$
cannot contract to the 22-dimensional orbit labelled $A_1$; therefore,
it maps to $2A_1$.

We deduce that the fibers above the $A_1$-orbit are 8-dimensional; by
semi-continuity again, the orbits labelled $3A_{1,b},3A_{1,c}$ map to
the orbit $3A_1$.

\lpara

Since by lemma \ref{orbites} the map $U$ of notation
\ref{U} is defined only on the open
orbit of $Hom(L_\alpha, \wedge^2 C_\alpha^*)$, the Mukai
flop is also defined only on the open orbit of $T^* Y$.
\qed

\end{document}